\documentclass[a4paper,12pt]{article} 

\usepackage{color}

\usepackage{amsfonts, amsmath, amsthm, amssymb}
\usepackage[T1]{fontenc}
\usepackage[cp1250]{inputenc}
\usepackage{xcolor}
\usepackage{graphicx}
\usepackage{amssymb}
\usepackage{amsmath}
\usepackage{mathptmx}
\usepackage{helvet}
\usepackage{courier}
\usepackage{txfonts}
\usepackage{tikz} 
\usetikzlibrary{arrows}
\usepackage{type1cm}

\usepackage{verbatim}

\usepackage{graphicx}
\usepackage{epsfig,amscd,amssymb,amsxtra,amsmath,amsthm}
\usepackage{type1cm}
\usepackage[T1]{fontenc}
\usepackage{graphics}
\usepackage[mathscr]{eucal}
\usepackage[all]{xy}
\usepackage{amsmath,amscd}

\newtheorem{theorem}{Theorem}[section]
\newtheorem{proposition}[theorem]{Proposition}
\newtheorem{definition}[theorem]{Definition}

\newtheorem{remark}[theorem]{Remark}

\newtheorem{corollary}[theorem]{Corollary}
\newtheorem{problem}[theorem]{Problem}
\newtheorem{observation}[theorem]{Observation}

\newcommand{\diam} {\mathop{\rm diam}\nolimits}

\newcommand{\Cl}  {\mathop{\rm Cl}\nolimits}

\begin{document}

\def\joinrel{\mkern-3mu}
\newcommand{\varproj}{\displaystyle \lim_{\multimapinv\joinrel-\joinrel-}}

\title{An uncountable family of smooth fans  that admit transitive homeomorphisms}
\author{Iztok Bani\v c, Goran Erceg, Judy Kennedy,  Chris Mouron and Van Nall}
\date{}

\maketitle

\begin{abstract}
We construct a family of uncountably many pairwise non-homeomorphic smooth fans that admit transitive homeomorphisms. In addition, we show that the star of Cantor fans admits a transitive homeomorphism.
\end{abstract}
\-

\section{Introduction}
Many examples of continua that admit transitive homeomorphisms may be found in the literature.  Most of the known examples of such continua have a very complicated topological structure, i.e., they are indecomposable or they are {decomposable} but have some other complicated property. However, smooth fans form a family of continua that have been considered not to be very complicated.  In previous papers, we show that the Cantor fan and the Lelek fan admit a transitive homeomorphisms, see \cite{BE,banic2} where more references may be found. We began to wonder if the Cantor fan and the Lelek fan were special in this regard among smooth fans. Gradually, we discovered more smooth fans that admit transitive homeomorphisms. We present in this paper a family of uncountably many pairwise non-homeomorphic smooth fans that admit transitive homeomorphisms.

We proceed as follows. In Section \ref{s1}, we introduce the definitions, notation and the well-known results that will be used later in the paper. In Section \ref{s2}, we show that the star of Cantor fans is another example of a smooth fan that admits a transitive homeomorphism, and then, in Section \ref{s3}, a family of uncountably many pairwise non-homeomorphic smooth fans that admit a transitive homeomorphism is constructed.

\section{Definitions and Notation}\label{s1}
The following definitions, notation and well-known results are needed in the paper.

\begin{definition}
Let $X$ and $Y$ be metric spaces, and let $f:X\rightarrow Y$ be a function.  We use  $\Gamma(f)=\{(x,y)\in X\times Y \ | \ y=f(x)\}$
to denote \emph{  the graph of the function $f$}.
\end{definition}

\begin{definition}
Let $X$ be a metric space, $x\in X$ and $\varepsilon>0$. We use $B(x,\varepsilon)$ to denote the open ball,  {centered} at $x$ with radius $\varepsilon$.
\end{definition}
\begin{definition}
We use $\mathbb N$ to denote the set of positive integers and $\mathbb Z$ to denote the set of integers.  
\end{definition}
\begin{definition}
Let $(X,d)$ be a compact metric space. Then we define \emph{$2^X$} by 
$$
2^{X}=\{A\subseteq X \ | \ A \textup{ is a non-empty closed subset of } X\}.
$$
Let $\varepsilon >0$ and let $A\in 2^X$. Then we define  \emph{$N_d(\varepsilon,A)$} by 
$$
N_d(\varepsilon,A)=\bigcup_{a\in A}B(a,\varepsilon).
$$
Let $A,B\in 2^X$. The function \emph{$H_d:2^X\times 2^X\rightarrow \mathbb R$}, defined by
$$
H_d(A,B)=\inf\{\varepsilon>0 \ | \ A\subseteq N_d(\varepsilon,B), B\subseteq N_d(\varepsilon,A)\},
$$
is called \emph{the Hausdorff metric}. The Hausdorff metric is in fact a metric and the metric space $(2^X,H_d)$ is called \emph{the hyperspace of the space $(X,d)$}. 
\end{definition}
\begin{remark}
Let $(X,d)$ be a compact metric space, let $A$ be a non-empty closed subset of $X$,  and let $(A_n)$ be a sequence of non-empty closed subsets of $X$. When we say $\displaystyle A=\lim_{n\to \infty}A_n$ with respect to the {Hausdorff} metric, we mean $\displaystyle A=\lim_{n\to \infty}A_n$ in $(2^X,H_d)$. 
\end{remark}
\begin{definition}
 \emph{A continuum} is a non-empty compact connected metric space.  \emph{A subcontinuum} is a subspace of a continuum, which is itself a continuum.
 \end{definition}
 \begin{definition}
Let $X$ be a continuum. 
\begin{enumerate}
\item The continuum $X$ is \emph{ unicoherent}, if for any subcontinua $A$ and $B$ of $X$ such that $X=A\cup B$,  the compactum $A\cap B$ is connected. 
\item The continuum $X$ is \emph{hereditarily unicoherent } provided that each of its subcontinua is unicoherent.
\item The continuum $X$ is a \emph{dendroid}, if it is an arcwise connected, hereditarily unicoherent continuum.
\item Let $X$ be a continuum.  If $X$ is homeomorphic to $[0,1]$, then $X$ is \emph{ an arc}.   
\item A point $x$ in an arc $X$ is called \emph{an end-point of the arc  $X$}, if  there is a homeomorphism $\varphi:[0,1]\rightarrow X$ such that $\varphi(0)=x$.
\item Let $X$ be a dendroid.  A point $x\in X$ is called an \emph{end-point of the dendroid $X$}, if for  every arc $A$ in $X$ that contains $x$, $x$ is an end-point of $A$.  The set of all end-points of $X$ will be denoted by $E(X)$. 
\item A continuum $X$ is \emph{a simple triod}, if it is homeomorphic to $([-1,1]\times {0})\cup (\{0\}\times [0,1])$.
\item A point $x$ in a simple triod $X$ is called \emph{the top-point} or just the \emph{top of the simple triod $X$}, if  there is a homeomorphism $\varphi:([-1,1]\times {0})\cup (\{0\}\times [0,1])\rightarrow X$ such that $\varphi(0,0)=x$.
\item Let $X$ be a dendroid.  A point $x\in X$ is called \emph{a ramification-point of the dendroid $X$}, if there is a simple triod $T$ in $X$ with the top   $x$.  The set of all ramification-points of $X$ will be denoted by $R(X)$. 
\item The continuum $X$ is \emph{a  fan}, if it is a dendroid with at most one ramification point $v$, which is called the top of the fan $X$ (if it exists).
\item Let $X$ be a fan.   For all points $x$ and $y$ in $X$, we define  \emph{$A[x,y]$} to be the arc in $X$ with end-points $x$ and $y$, if $x\neq y$. If $x=y$, then we define $A[x,y]=\{x\}$.
\item Let $X$ be a fan with the top $v$. We say that that the fan $X$ is \emph{smooth} if for any $x\in X$ and for any sequence $(x_n)$ of points in $X$,
$$
\lim_{n\to \infty}x_n=x \Longrightarrow \lim_{n\to \infty}A[v,x_n]=A[v,x]
$$ 
with respect to the Hausdorff metric.
\item Let $X$ be a fan.  We say that $X$ is \emph{a Cantor fan}, if $X$ is homeomorphic to the continuum
$$
\bigcup_{c\in C}{S}_c,
$$
where $C\subseteq [0,1]$ is the standard Cantor set and for each $c\in C$, ${S}_c$ is the straight line segment in the plane from $(0,0)$ to $(c,1)$. See Figure \ref{fig000}, where a Cantor fan is pictured.
\begin{figure}[h!]
	\centering
		\includegraphics[width=25em]{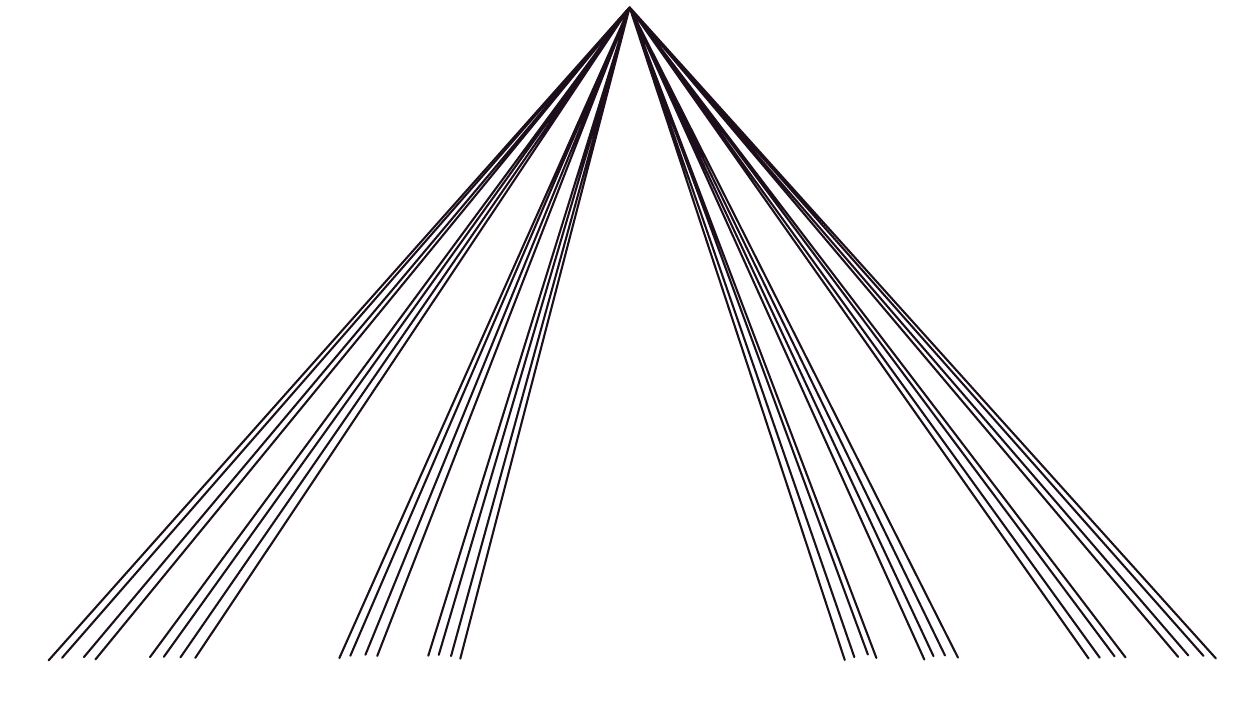}
	\caption{A Cantor fan}
	\label{fig000}
\end{figure}  
\item Let $X$ be a fan.  We say that $X$ is \emph{a Lelek fan}, if it is smooth and $\Cl(E(X))=X$. See Figure \ref{figure2}, where a Lelek fan is pictured.
\begin{figure}[h!]
	\centering
		\includegraphics[width=25em]{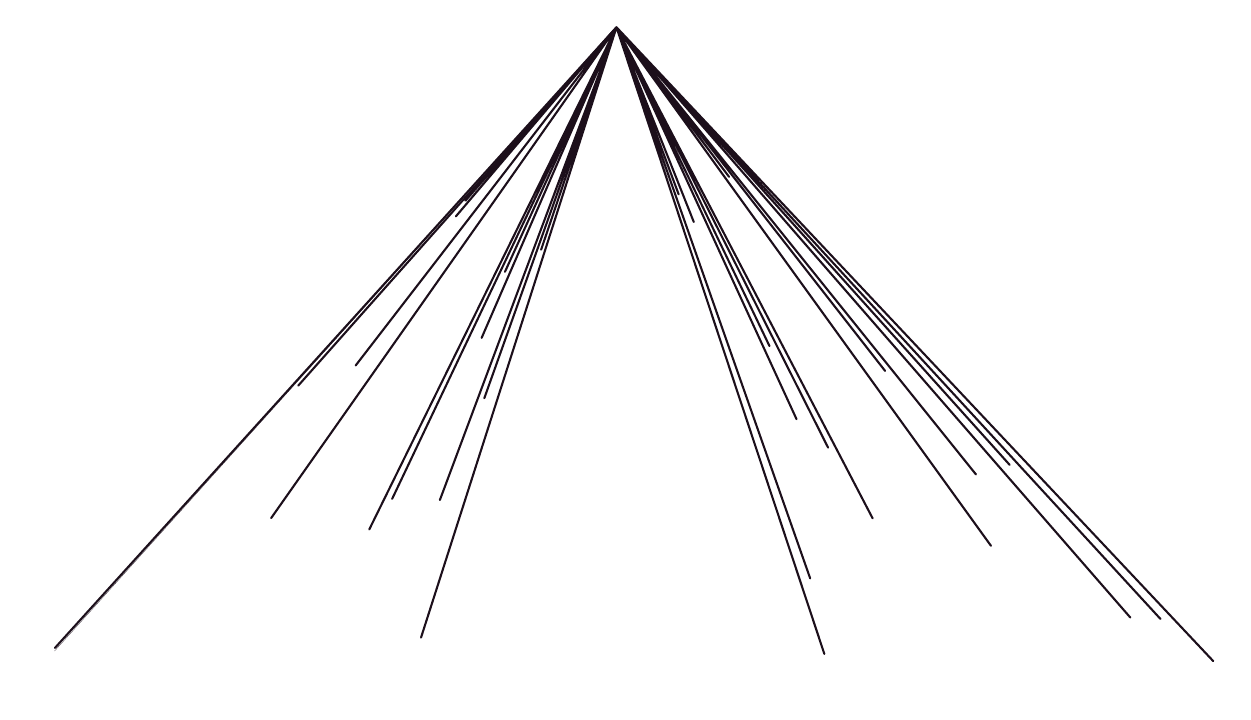}
	\caption{A Lelek fan}
	\label{figure2}
\end{figure} 
\end{enumerate}
\end{definition}
\begin{observation}
	It is a well-known fact that the Cantor fan is smooth and that any subcontinuum of a smooth fan is itself a smooth fan. 
\end{observation}
An example of a Lelek fan was constructed by A.~Lelek in \cite{lelek}.  He also showed that the set of the end-points of any Lelek fan is a dense one-dimensional set in it. Also, it is the only non-degenerate smooth fan with a dense set of end-points.  This was proved independently by W.~D.~Bula and L.~Oversteegen  in \cite{oversteegen} and by W. ~Charatonik in  \cite{charatonik}.  See \cite{nadler} for more information about continua, fans and their properties.
\begin{definition}
	Let $X$ and $Y$ be any continua and let $f:X\rightarrow Y$ be a continuous mapping. We say that $f$ is \emph{confluent}, if for every subcontinuum $S$ of $Y$ and for each component $C$ of $f^{-1}(S)$,
	$$
	f(C)=S.
	$$
\end{definition}
The following is a well-known result.
\begin{theorem}\label{char}
	Let $X$ and $Y$ be any continua and let $f:X\rightarrow Y$ be a confluent surjection. If $X$ is a smooth fan, then also $Y$ is a smooth fan.
	\end{theorem}
\begin{proof}
	See \cite[Theorem 13, page 33]{charatonik3}.
\end{proof}
\begin{definition}
Let $(X,f)$ be a dynamical system.  We say that $(X,f)$ is  \emph{transitive}, if for all non-empty open sets $U$ and $V$ in $X$,  there is a non-negative integer $n$ such that $f^n(U)\cap V\neq \emptyset$. We say that the mapping $f$ is \emph{transitive}, if $(X,f)$ is transitive.
\end{definition}
\begin{definition}\label{povezava}
Let $X$ be a compact metric space. We say that $X$  \emph{admits a transitive homeomorphism}, if there is a homeomorphism $f:X\rightarrow X$ such that $(X,f)$ is transitive.
\end{definition}
\begin{definition}
	For non-empty compact metric spaces $X$ and $Y$, we use 
$p_1:X\times Y\rightarrow X$ and $p_2:X\times Y\rightarrow Y$
to denote \emph{the standard projections} defined by
$p_1(s,t)=s$ and $p_2(s,t)=t$ for all $(s,t)\in X\times Y$.
\end{definition}

\section{A star of Cantor fans}\label{s2}
In this section, we construct an example of a smooth fan, the star of Cantor fans, and show that it admits a transitive homeomorphism.   
\begin{definition}\label{mu}
Let $F$ be a smooth fan and let $(F_n)$ be a sequence of smooth fans in the plane such that 
\begin{enumerate}
    \item for each positive integer $n$, $F_n$ is homeomorphic to $F$,
	\item for each positive integer $n$, $\diam(F_n)\leq \frac{1}{2^n}$,
	\item for each positive integer $n$, $(0,0)$ is the top of $F_n$, and 
	\item for all positive integers $m$ and $n$, $F_m\cap F_n=\{(0,0)\}$. 
\end{enumerate}
Any space $X$ that is homeomorphic to $\bigcup_{n=1}^{\infty}F_n$, is called a star of $F$'s. 
\end{definition}
\begin{observation}
Let $F$ be a smooth fan and let $X$ be a star of $F$'s. Then $X$ is also a smooth fan.
\end{observation}
\begin{definition}
Let $F$ be a smooth fan and let $X$ be a star of $F$'s. If $F$ is a
\begin{enumerate}
	\item Cantor fan, then $X$ is called a star of Cantor fans, see Figure \ref{figure3}.
	\begin{figure}[h!]
	\centering
		\includegraphics[width=30em]{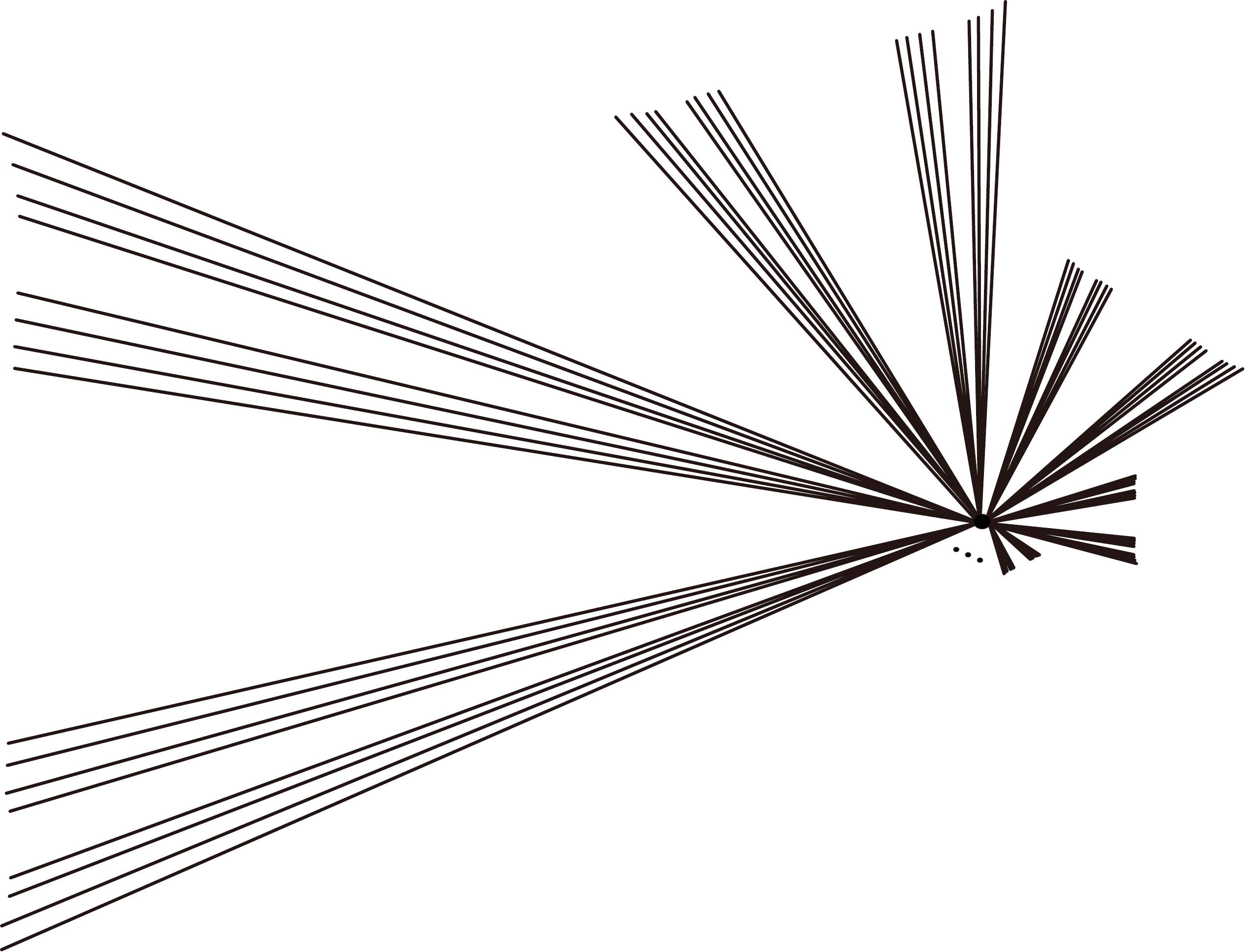}
	\caption{The star of Cantor fans}
	\label{figure3}
\end{figure}  
	\item Lelek fan, then $X$ is called a star of Lelek fans.
	\item star of Cantor fans, then $X$ is called a star of stars of Cantor fans.
\end{enumerate}
\end{definition}
\begin{observation}\label{osje}
Note that any star of Lelek fans is again a Lelek fan and that any star of stars of Cantor fans is again a star of Cantor fans.	
Also, note that any two stars of Cantor fans are homeomorphic. 
\end{observation}

\begin{definition}
	We use
	\begin{enumerate}
		\item $C$ to denote the standard middle-third Cantor set in $[0,1]$, 
		\item $I$ to denote the closed interval $[0,1]$,
		\item $P$ to denote the topological product $P=C\times I$, 
		\item $R$ to denote the subspace $R=P\cap \{(s,t)\in I\times I \ | \ s\leq t\}$
		of the space $P$, and
		\item $\varphi$ to denote the function $\varphi:P\rightarrow R$ that is defined by
		$
		\varphi(c,t)=(c,c\cdot t)
		$
		for each $(c,t)\in P$.
	\end{enumerate} 
	see Figure \ref{figure4}.
	\begin{figure}[h!]
	\centering
		\includegraphics[width=30em]{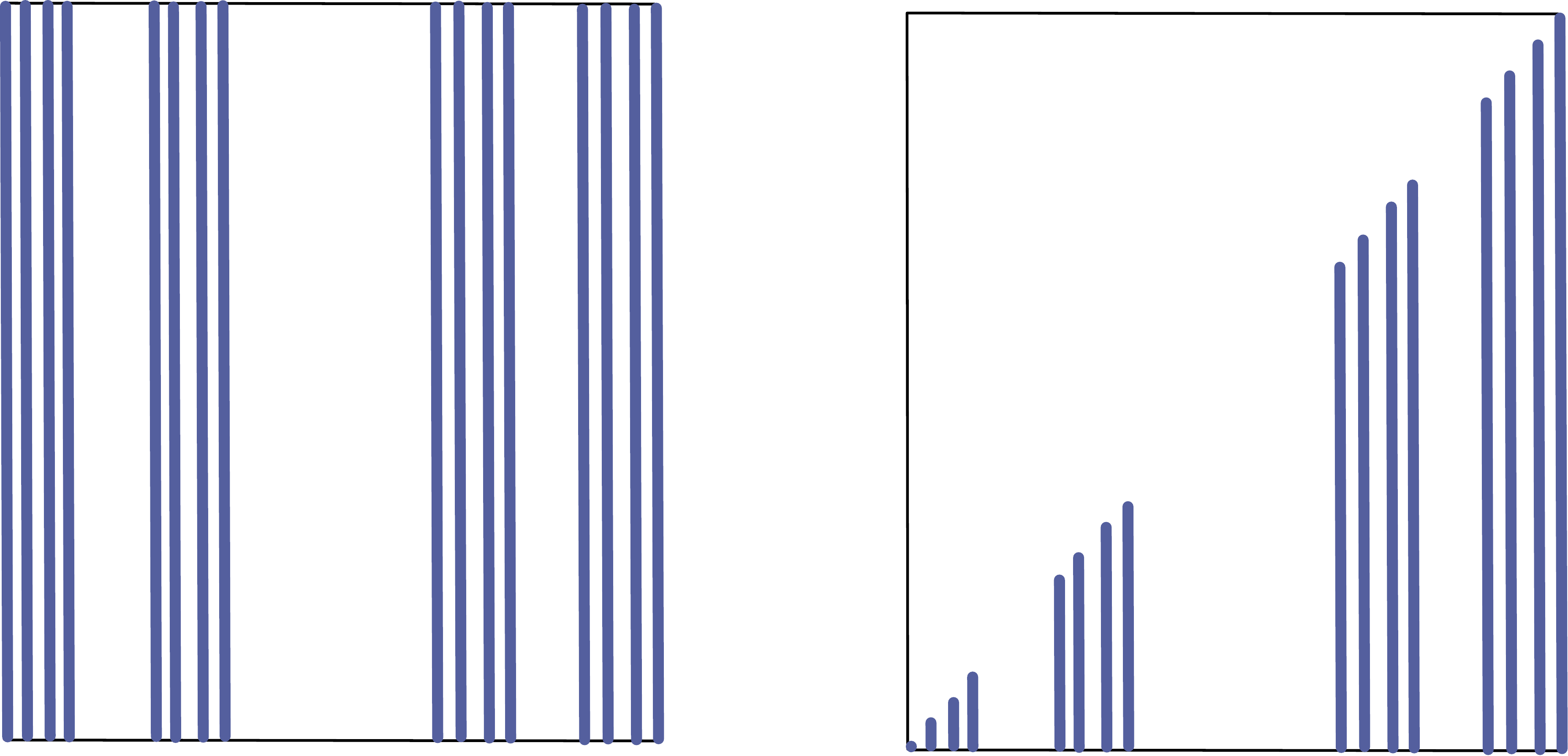}
	\caption{The spaces $P$ and $R$}
	\label{figure4}
\end{figure}  
\end{definition}
\begin{observation}
	 Note that $\varphi$ is a continuous surjection. 
\end{observation}

\begin{definition}
	For any function $f:P\rightarrow P$, we denote by $f_R$ the function $f_R:R\rightarrow R$, defined by
	$$
	f_R(c,t)=\begin{cases}
				(0,0)\text{;} & c=0 \\
				\varphi(f(\varphi^{-1}(c,t)))\text{;} & c\neq 0
			\end{cases}
	$$
	for any $(c,t)\in R$. 
\end{definition}
\begin{proposition}\label{hlavna}
	Let $f:P\rightarrow P$ be a function such that $f(\{0\}\times I)\subseteq \{0\}\times I$ and $f((C\setminus \{0\})\times I)\subseteq (C\setminus \{0\})\times I$. Then the following hold.
	\begin{enumerate}
		\item If $f$ is surjective, then $f_R$ is surjective.
	 	\item If $f$ is injective, then $f_R$ is injective.
		\item If $f$ is continuous, then $f_R$ is continuous.
	 	\item If $f$ is a homeomorphism, then $f_R$ is a homeomorphism.
	 	\item If $f$ is transitive, then $f_R$ is transitive.
	\end{enumerate}
\end{proposition}
\begin{proof}
First, suppose that $f$ is surjective. To show that $f_R$ is surjective, let $(c,t)\in R$. Also, let $(c_0,t_0)\in \varphi^{-1}(c,t)$. Since $f$ is surjective, there is a point $(c_1,t_1)\in P$ such that $f(c_1,t_1)=(c_0,t_0)$.  We treat the following possible cases.
\begin{enumerate}
	\item $c_1\neq 0$. It follows that 
$$
f_R(\varphi(c_1,t_1))=\varphi(f(\varphi^{-1}(\varphi(c_1, t_1))))=\varphi(f(c_1,t_1))=\varphi(c_0,t_0)=(c,t).
$$
	\item $c_1=0$. Then $c_0=0$ and it follows that $(c,t)=(0,0)$. Therefore, 
	$$
	f_R(\varphi(c_1,t_1))=f_R(\varphi(0,t_1))=f_R(0,0)=(0,0)=(c,t).
	$$
\end{enumerate}
It follows that $f_R$ is surjective.

Next, suppose that $f$ is injective. To see that $f_R$ is injective, let $(c_1,t_1),(c_2,t_2)\in R$ be such points that $f_R(c_1,t_1)=f_R(c_2,t_2)$. To see that $(c_1,t_1)=(c_2,t_2)$, we consider the following possible cases.
\begin{enumerate}
	\item $c_1=0$. It follows that $t_1=0$ and
	$$
	f_R(c_1,t_1)=f_R(0,0)=(0,0).
	$$
	Therefore, $f_R(c_2,t_2)=(0,0)$. Suppose that $c_2\neq 0$. Then 
	$$
	f_R(c_2,t_2)=\varphi(f(\varphi^{-1}(c_2,t_2)))=\varphi\Big(f\Big(c_2,\frac{t_2}{c_2}\Big)\Big).
	$$
	  Since $c_2\neq 0$ and $f((C\setminus \{0\})\times I)\subseteq (C\setminus \{0\})\times I$, it follows that $\varphi(f(c_2,\frac{t_2}{c_2}))\neq (0,0)$. Therefore, $f_R(c_1,t_1)\neq (0,0)$, which is a contradiction. Therefore, $c_2=0$ and it follows that also $t_2=0$. Hence, $(c_1,t_1)=(c_2,t_2)$.
	\item $c_1\neq 0$. If $c_2=0$, we obtain a contradiction similarly as in the previous case. Therefore, $c_2\neq 0$. It follows from $f_R(c_1,t_1)=f_R(c_2,t_2)$ that $\varphi(f(\varphi^{-1}(c_1,t_1)))=\varphi(f(\varphi^{-1}(c_2,t_2)))$. Therefore, $\varphi(f(c_1,\frac{t_1}{c_1}))=\varphi(f(c_2,\frac{t_2}{c_2}))$. Since $f$ is injective,  since $f((C\setminus \{0\})\times I)\subseteq (C\setminus \{0\})\times I$, and since $\varphi$ restricted to $(C\setminus \{0\})\times I$ is injective, it follows that $(c_1,\frac{t_1}{c_1})=(c_2,\frac{t_2}{c_2})$. Therefore, $(c_1,t_1)=(c_2,t_2)$.
\end{enumerate}
Thus $f_R$ is injective.

	Next, suppose that $f$ is continuous and let $(c_n,t_n)$ be a sequence of points in $R$ such that $\displaystyle	\lim_{n\to \infty}(c_n,t_n)=(0,0)$ and such that for each positive integer $n$, $(c_n,t_n)\neq (0,0)$. 
	Then 
	\begin{align*}
		&\lim_{n\to \infty}f_R(c_n,t_n)=\lim_{n\to \infty}\varphi(f(\varphi^{-1}(c_n,t_n)))=\lim_{n\to \infty}\varphi\left(f\Big(c_n,\frac{t_n}{c_n}\Big)\right)=
		\\
		&\lim_{n\to \infty}\varphi\left(p_1\Big(f\Big(c_n,\frac{t_n}{c_n}\Big)\Big),p_2\Big(f\Big(c_n,\frac{t_n}{c_n}\Big)\Big)\right)=
		\\
		&\lim_{n\to \infty}\left(p_1\Big(f\Big(c_n,\frac{t_n}{c_n}\Big)\Big),p_1\Big(f\Big(c_n,\frac{t_n}{c_n}\Big)\Big)\cdot p_2\Big(f\Big(c_n,\frac{t_n}{c_n}\Big)\Big)\right)=(0,0)
	\end{align*}
	since $\displaystyle \lim_{n\to \infty}p_1\Big(f\Big(c_n,\frac{t_n}{c_n}\Big)\Big)=0$ and since the sequence $\left(p_2\Big(f\Big(c_n,\frac{t_n}{c_n}\Big)\Big)\right)$ is bounded (by $0$ from below and by $1$ from above). It follows that $f_R$ is continuous. 
	
	Next, suppose that $f$ is a homeomorphism. It follows from the previous claims that also $f_R$ is a homeomorphism.  
	
	Finally, suppose that $f$ is transitive. Let $U$ and $V$ be non-empty open sets in $R$ and let $U'=U\setminus \{(0,0)\}$ and $V'=V\setminus \{(0,0)\}$. Since $R$ does not have any isolated points, it follows that $U'$ and $V'$ are also non-empty open sets in $R$. Since $\varphi$ is continuous, it follows that $\varphi^{-1}(U')$ and $\varphi^{-1}(V')$ are open in $P$ and it follows from the definition of $\varphi$ that $\varphi^{-1}(U')\cap (\{0\}\times I)=\emptyset$ and $\varphi^{-1}(V')\cap (\{0\}\times I)=\emptyset$. Since $f$ is transitive, there is a non-negative integer $n$ such that $f^n(\varphi^{-1}(U'))\cap \varphi^{-1}(V')\neq \emptyset$. Let $n$ be such a non-negative integer and let $(c,t)\in f^n(\varphi^{-1}(U'))\cap \varphi^{-1}(V')$. Since $(c,t)\in \varphi^{-1}(V')$, it follows that $c\neq 0$. It follows from $\varphi(c,t)\in V'$ and $\varphi(c,t)=(c,c\cdot t)$ that $(c,c\cdot t)\in V'$. Then, using $f((C\setminus \{0\})\times I)\subseteq (C\setminus \{0\})\times I$, we get that 
\begin{align*}
	&f_R^n(c,c\cdot t)=(\varphi\circ f \circ \varphi^{-1})^n(c,c \cdot t)=(\varphi\circ f^n \circ \varphi^{-1})(c,c \cdot t)=
	\\
	&(\varphi\circ f^n)(c,t)\in \varphi(\varphi^{-1}(U'))\subseteq U'.
\end{align*}
Therefore, $(c,c\cdot t)\in f_R^n(U')\cap V'$. Since $f_R^n(U')\cap V'\subseteq f_R^n(U)\cap V$, it follows that $f_R^n(U)\cap V\neq \emptyset$. Therefore, $f_R$ is transitive.
\end{proof}
\begin{observation}\label{lipa}
	Note that there is a transitive homeomorphism $f:P\rightarrow P$ such that $f(\{0\}\times I)\subseteq \{0\}\times I$ and $f((C\setminus \{0\})\times I)\subseteq (C\setminus \{0\})\times I$. One such homeomorphism can be constructed using \cite[Theorem 3.32, page 17]{BE}, where a topological conjugacy of such a homeomorphism is obtained.
\end{observation}
\begin{definition}
	We use
	\begin{enumerate}
		\item $\sim$ to denote the relation $\sim$ on $P$, which is defined by
		$$
		(c_1,t_1)\sim(c_2,t_2) ~~~ \Longleftrightarrow ~~~ (c_1,t_1)=(c_2,t_2) \textup{ or } t_1=t_2=0 
		$$
		for all $(c_1,t_1),(c_2,t_2)\in P$. 
		\item $q$ to denote the quotient function $q:P\rightarrow P/_{\sim}$, defined by
		$$
		q(c,t)=[(c,t)]
		$$
		for each $(c,t)\in P$.
		\item $\sim_{R}$ to denote the restriction of the relation $\sim$ to $R$.
	\end{enumerate} 
\end{definition}
\begin{observation}
	 Note that $P/_{\sim}$ is a Cantor fan and that $R/_{\sim_R}$ is a star of Cantor fans. 
\end{observation}
\begin{definition}
Let $X$ be a compact metric space, let $\sim$ be an equivalence relation on $X$,  and let $f:X\rightarrow X$ be a function such that for all $x,y\in X$,
$$
x\sim y  \Longleftrightarrow f(x)\sim f(y).
$$
 Then we let $f^{\star}:X/_{\sim}\rightarrow X/_{\sim}$ be defined by   
$
f^{\star}([x])=[f(x)]
$
for any $x\in X$. 
\end{definition}
The following proposition is a well-known result.
\begin{proposition}\label{kvocienti}
Let $X$ be a compact metric space, let $\sim$ be an equivalence relation on $X$, and  let $f:X\rightarrow X$ be a function such that for all $x,y\in X$,
$$
x\sim y  \Longleftrightarrow f(x)\sim f(y).
$$
 Then the following hold.
\begin{enumerate}
\item $f^{\star}$ is a well-defined function from  $X/_{\sim}$ to $X/_{\sim}$. 
\item If $f$ is continuous, then $f^{\star}$ is continuous.
\item If $f$ is a homeomorphism, then $f^{\star}$ is a homeomorphism.
\item If $f$ is transitive, then $f^{\star}$ is transitive.
\end{enumerate}
\end{proposition}
\begin{proof}
	See \cite[Theorem 3.4]{BE}. 
\end{proof}
\begin{theorem}\label{starcek}
	The star of Cantor fans admits a transitive homeomorphism.
\end{theorem}
\begin{proof}
	By Observation \ref{lipa}, there is a transitive homeomorphism $f:P\rightarrow P$ such that 
	$$
	f(\{0\}\times I)\subseteq \{0\}\times I
	$$ 
	and 
	$$f((C\setminus \{0\})\times I)\subseteq (C\setminus \{0\})\times I.
	$$
	 Fix such a homeomorphism $f$. By Proposition \ref{hlavna}, $f_R$ is a transitive homeomorphism from $R$ to $R$. It follows from Proposition \ref{kvocienti} that $f_R^{\star}$ is a transitive homeomorphism from $R/_{\sim_{R}}$ to $R/_{\sim_{R}}$.
\end{proof}
Suppose that $F$ is either a Cantor fan, a Lelek fan or a star of Cantor fans. Then, as seen above, a star of $F$'s also admits a transitive homeomorphism (by Theorem \ref{starcek} and Observation \ref{osje}, since both the Cantor fan and the Lelek fan admit a transitive homeomorphism as seen in \cite{BE} and in \cite{banic2}). Therefore,  the following open problems are a good place to finish this section.  
\begin{problem}
	Let $F$ be a smooth fan that admits a transitive homeomorphism. Does the star of $F$'s also admit a transitive homeomorphism?
\end{problem}
\begin{problem}
	Let $F$ be a smooth fan. If a star of $F$'s admits a transitive homeomorphism, then does $F$ admit a transitive homeomorphism?
\end{problem}

\section{An uncountable family of smooth fans that admit transitive homeomorphisms}\label{s3}
In this section, we present our main result of the paper, an uncountable family of pairwise non-homeomorphic smooth fans that admit transitive homeomorphisms. First, we construct, using a Mahavier product of a closed relation, a space that is homeomorphic to a subspace of the product of a Cantor set and an interval in such a way that the shift map on it is a transitive homeomorphism. Second, with some identifications, we get uncountably many smooth fans while keeping the transitivity of the induced homeomorphisms.  Great care must be taken to see that the model for the Mahavier product that we present is what we claim it is.  Identifications must also be done with care to ensure the transitivity of the induced homeomorphisms. 

We begin with the following definitions.
\begin{definition}
Let $X$ be a {non-empty} compact metric space and let ${F}\subseteq X\times X$ be a relation on $X$. If ${F}$ is closed in $X\times X$, then we say that ${F}$ is  \emph{  a closed relation on $X$}.  
\end{definition}
\begin{definition}
Let $X$ be a {non-empty} compact metric space and let ${F}$ be a closed relation on $X$. For each positive integer $m$, we call 
$$
X_F^m=\Big\{(x_1,x_2,\ldots ,x_{m+1})\in \prod_{{ i={1}}}^{{m{ +1}}}X \ | \ \textup{ for each } i\in{  \{{1,2},\ldots ,m\}}, (x_{i},x_{i+1})\in {F}\Big\}
$$
 \emph{ the $m$-th Mahavier product of ${F}$}, and we call
$$
X_F^+=\Big\{(x_1,x_2,x_3,\ldots )\in \prod_{{ i={1}}}^{\infty}X \ | \ \textup{ for each { positive} integer } i, (x_{i},x_{i+1})\in {F}\Big\}
$$
\emph{ the  Mahavier product of ${F}$}, and 
$$
X_F=\Big\{(\ldots,x_{-3},x_{-2},x_{-1},{x_0}{ ;}x_1,x_2,x_3,\ldots )\in \prod_{i={-\infty}}^{\infty}X \ | \ \textup{ for each  integer } i, (x_{i},x_{i+1})\in {F}\Big\}
$$
\emph{ the two-sided  Mahavier product of ${F}$}.
\end{definition}

\begin{definition}
Let $X$ be a {non-empty} compact metric space and let ${F}$ be a closed relation on $X$. 
The function  $\sigma_F^{+} : {X_F^+} \rightarrow {X_F^+}$, 
 defined by 
$$
\sigma_F^{+} ({x_1,x_2,x_3,x_4},\ldots)=({x_2,x_3,x_4},\ldots)
$$
for each $({x_1,x_2,x_3,x_4},\ldots)\in {X_F^+}$, 
is called \emph{   the shift map on ${X_F^+}$}.      
The function  $\sigma_F : {X_F} \rightarrow {X_F}$, 
 defined by 
$$
\sigma_F (\ldots,x_{-3},x_{-2},x_{-1},{x_0};x_1,x_2,x_3,\ldots )=(\ldots,x_{-3},x_{-2},x_{-1},{x_0},x_1;x_2,x_3,\ldots )
$$
for each $(\ldots,x_{-3},x_{-2},x_{-1},{x_0};x_1,x_2,x_3,\ldots )\in {X_F}$, 
is called \emph{   the shift map on ${X_F}$}.    
\end{definition}
\begin{observation}\label{juju}
Note that $\sigma_F$ is always a homeomorphism while $\sigma_F^+$ may not be a homeomorphism.
\end{observation}
Next, we define a space $\mathbb X$, which will be used to construct our uncountable family of pairwise non-homeomorphic smooth fans that admit transitive homeomorphisms.  
\begin{definition}
	We use $\mathbb X$ to denote the set
	$$
	\mathbb X=([0,1]\cup[2,3]\cup[4,5]\cup[6,7]\cup\ldots)\cup\{\infty\}.
	$$
	We equip $\mathbb X$ with the Alexandroff one-point compactification topology $\mathcal T$; i.e., $\mathcal T$ is obtained on $\mathbb X$ from the Alexandroff one-point compactification (also known as the Alexandroff extension) of the space $[0,1]\cup[2,3]\cup[4,5]\cup[6,7]\cup\ldots$ (which is a subspace of the Euclidean line $\mathbb R$) with the point $\infty$. See \cite[pages 166--171]{engelking1} or \cite[pages 135--145]{willard} for more information on (one-point) compactifications. 
	\end{definition}
	\begin{observation}\label{prej}
	For each non-negative integer $k$, let $q_k=1-\frac{1}{2^k}$ and let 
	$$
	X=[q_0,q_1]\cup [q_2,q_3]\cup [q_4,q_5]\cup [q_6,q_7]\cup \ldots \{1\}
	$$ 
	(we equip $X$ with the usual topology). Note that the compacta $\mathbb X$ and $X$ are homeomorphic. 
	\end{observation}
	\begin{definition}
		Let $X$ be the compactum from Observation \ref{prej} and let $h:X\rightarrow \mathbb X$ be any homeomorphism such that for each non-negative integer $k$, $h(q_k)=k$. On the space $\mathbb X$, we always use the metric $\textup{d}_{\mathbb X}$ that is defined by 
	$$
	\textup{d}_{\mathbb X}(x,y)=|h^{-1}(y)-h^{-1}(x)|
	$$
	 for all $x,y\in \mathbb X$.
	\end{definition}
	\begin{observation}
	Note that the topology $\mathcal T_{\textup{d}_{\mathbb X}}$ on $\mathbb X$, that is obtained from the metric $\textup{d}_{\mathbb X}$, is exactly the one-point compactification topology $\mathcal T$ on $\mathbb X$. 	Also, note that (in this setting) for each non-negative integer $k$, 
	$$	
	\diam([2k,2k+1])=\frac{1}{2^{2k+1}}.
	$$
	\end{observation}
	
	\begin{definition}
		For each non-negative integer $k$, we use $I_{k+1}$ to denote  
		$$
			I_{k+1}=[2k,2k+1].
			$$
	\end{definition}
	\begin{observation}
		Note that for each positive integer $k$, 
		$$
		\diam(I_k)=\frac{1}{2^{2k-1}}.
		$$
	\end{observation}
	\begin{definition}
		We use the product metric $\textup{D}_{\mathbb X}$ on the  product $\prod_{k=-\infty}^{\infty}\mathbb X$, which is defined by 
		$$
		\textup{D}_{\mathbb X}(\mathbf x,\mathbf y)=\sup \Big\{\frac{\textup{d}_{\mathbb X}(\mathbf x(k),\mathbf y(k))}{2^{|k|}} \ \big| \ k \textup{ is an integer}\Big\}
		$$
		for all $\mathbf x,\mathbf y\in \prod_{k=-\infty}^{\infty}\mathbb X$. 
	\end{definition}
	\begin{observation}
		Since $\mathbb X$ is compact it follows that for all $\mathbf x,\mathbf y\in \prod_{k=-\infty}^{\infty}\mathbb X$, 
		$$
		\sup \Big\{\frac{\textup{d}_{\mathbb X}(\mathbf x(k),\mathbf y(k))}{2^{|k|}} \ \big| \ k \textup{ is an integer}\Big\}=\max \Big\{\frac{\textup{d}_{\mathbb X}(\mathbf x(k),\mathbf y(k))}{2^{|k|}} \ \big| \ k \textup{ is an integer}\Big\}
		$$
		and, therefore, for all $\mathbf x,\mathbf y\in \prod_{k=-\infty}^{\infty}\mathbb X$, 

		$$
		\textup{D}_{\mathbb X}(\mathbf x,\mathbf y)=\max \Big\{\frac{\textup{d}_{\mathbb X}(\mathbf x(k),\mathbf y(k))}{2^{|k|}} \ \big| \ k \textup{ is an integer}\Big\}.
		$$
	\end{observation}
	Next, we define the closed relation $H$ on $\mathbb X$ that will play an importaint role in our construction of an uncountable family of pairwise non-homeomorphic smooth fans that admit transitive homeomorphisms. 
	\begin{definition} 
	We use $H$ to denote the closed relation on $\mathbb X$ that is defined as follows:
	\begin{align*}
		H=&\Big\{\big(t,t^{\frac{1}{3}}\big) \ \big| \ t\in I_1\Big\}\cup \Big\{\big(t,(t-2)^2+2\big) \ \big| \ t\in I_2\Big\}\cup  \\
		&\Big\{(t,t+2) \ \big| \ t\in I_1\cup I_2\cup I_3\cup I_4\cup\ldots\Big\}\cup\\
		&\Big\{(t,t-2) \ \big| \ t\in I_2\cup I_3\cup I_4\cup I_5\cup\ldots\Big\}\cup\\
		&\Big\{(t,t) \ \big| \ t\in I_3\cup I_4\cup I_5\cup I_6\cup\ldots\Big\}\cup \Big\{(\infty,\infty)\Big\};
	\end{align*}
	see Figure \ref{uncun}. 
	\begin{figure}[h!]
	\centering
		\includegraphics[width=30em]{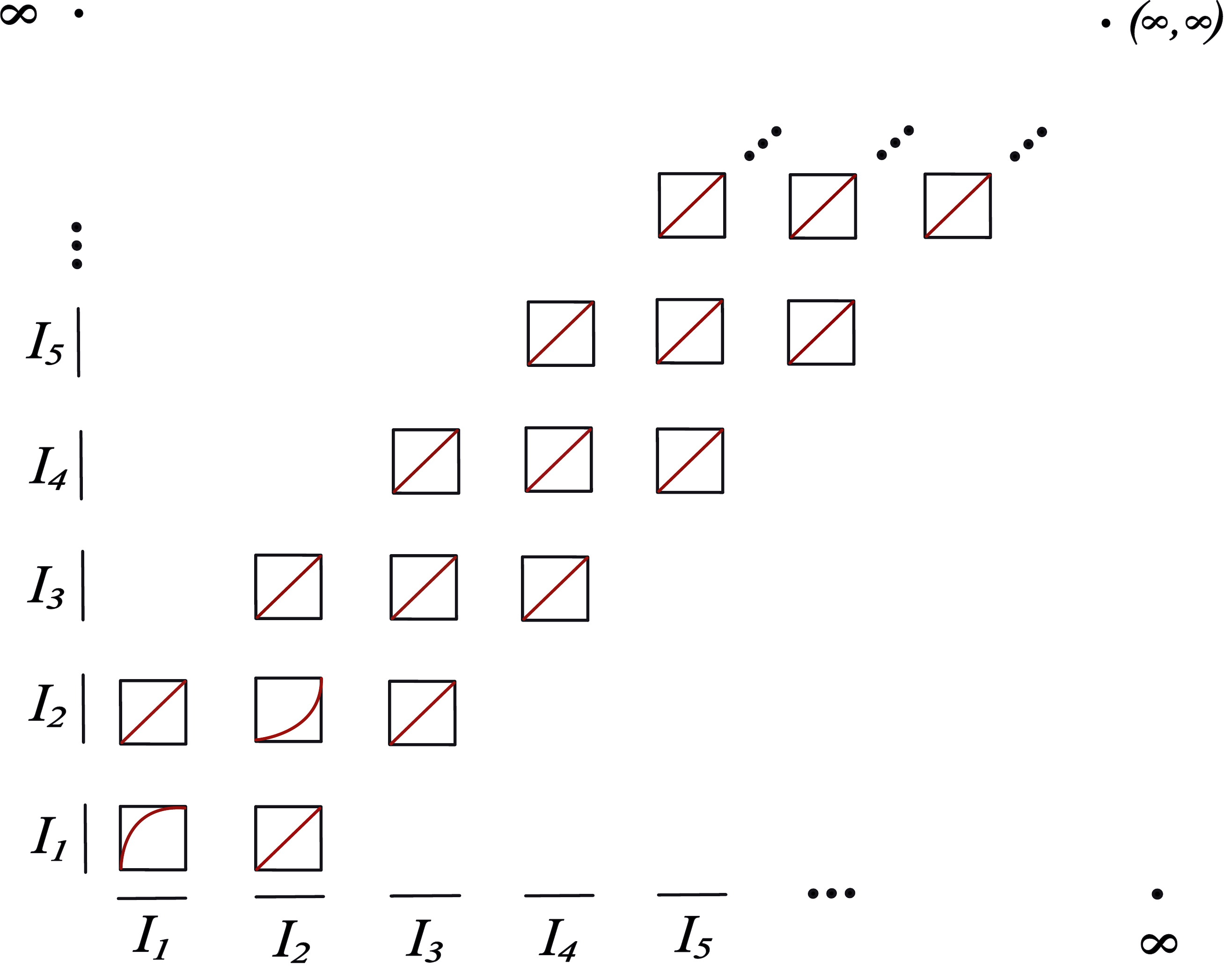}
	\caption{The relation $H$ on $\mathbb X$}
	\label{uncun}
\end{figure}
 We also use $\sigma_H^+$ to denote the shift map on the Mahavier product $\mathbb X_H^+$ and $\sigma_H$ to denote the shift map on the two-sided Mahavier product $\mathbb X_H$. 
\end{definition}
First, we prove that the shift map $\sigma_H$ is transitive (Theorem \ref{lejgaungatamlehre}). To do that, we use Theorems \ref{A} and \ref{B} from \cite{BE}.
\begin{definition}
Let $X$ be a compact metric space,  let $F$ be a closed relation on $X$ and let $x\in X$. Then we define 
$$
\mathcal U^{\oplus}_F(x)=\{y\in X \ | \ \textup{there are } n\in \mathbb N \textup{ and } \mathbf x\in X_F^{n} \textup{ such that } \mathbf x(1)=x, \mathbf x(n)=y \}
$$
and we call it the forward impression of $x$ by $F$.
\end{definition}
\begin{theorem}\label{A}
{ Let $X$ be a compact metric space, let $F$ be a closed relation on $X$,  let  $\{f_{\alpha} \ | \ \alpha \in A\}$ be a non-empty collection of continuous functions from $X$ to $X$ such that $F^{-1}=\bigcup_{\alpha\in A}\Gamma(f_{\alpha})$,  and let  $\{g_{\beta} \ | \ \beta \in B\}$ be a non-empty collection of continuous functions from $X$ to $X$ such that $F=\bigcup_{\beta\in B}\Gamma(g_{\beta})$. } If there is a dense set $D$ in $X$ such that for each $s\in D$,  $\Cl(\mathcal U^{\oplus}_F(s))=X$, then $(X_F^+,\sigma_F^+)$ is transitive. 
\end{theorem}
\begin{proof}
	See \cite[Theorem 4.8, page 18]{BE}.
\end{proof}
\begin{theorem}\label{B}
Let $X$ be a compact metric space and let $F$ be a closed relation on $X$ such that $p_1(F)=p_2(F)=X$. The following statements are equivalent. 
\begin{enumerate}
\item The map $\sigma_F^+$ is transitive.
\item The homeomorphism $\sigma_F$ is transitive. 
\end{enumerate}
\end{theorem}
\begin{proof}
	The theorem follows from \cite[Theorem 4.5, page 17]{BE}.
\end{proof}
\begin{theorem}\label{lejgaungatamlehre}
	The dynamical system $(\mathbb X_H,\sigma_H)$ is transitive.
\end{theorem}
\begin{proof}
	Since $p_1(H)=p_2(H)=\mathbb X$, it suffices to see that $(\mathbb X_H^+,\sigma_H^+)$ is transitive (by Theorem \ref{B}). We use Theorem \ref{A} to do so. Let $f_1,f_2,f_3:\mathbb X\rightarrow \mathbb X$ be defined as follows.
	For each $x\in \mathbb X$, let 
	$$
	f_1(x)=\begin{cases}
				x^{\frac{1}{3}}\text{;} & x\in [0,1] \\
				x^{2}\text{;} & x\in [2,3]  \\
				x\text{;} & x\in \mathbb X\setminus ([0,1]\cup[2,3]),
			\end{cases}
	$$
	$$
	f_2(x)=\begin{cases}
				x+2\text{;} & x\in \bigcup_{k=0}^\infty[4k,4k+1] \\
				x-2\text{;} & x\in \bigcup_{k=0}^\infty[4k+2,4k+3]  \\
				\infty\text{;} & x=\infty,
			\end{cases}
	$$
	and
	$$
	f_3(x)=\begin{cases}
	            x^{\frac{1}{3}}\text{;} & x\in [0,1] \\
				x+2\text{;} & x\in \bigcup_{k=0}^\infty[4k+2,4k+3] \\
				x-2\text{;} & x\in \bigcup_{k=0}^\infty[4k+4,4k+5]  \\
				\infty\text{;} & x=\infty. 
			\end{cases}
	$$
	Note that $f_1$, $f_2$, and $f_3$ are homeomorphisms from $\mathbb X$ to $\mathbb X$ such that $H=\Gamma(f_1) \cup \Gamma (f_2) \cup \Gamma (f_3)$. Similarly, $H^{-1}=\Gamma(f_1^{-1}) \cup \Gamma (f_2^{-1}) \cup \Gamma (f_3^{-1})$. So, all the initial conditions from Theorem \ref{A} are satisfied. To see that $(\mathbb X_H^+,\sigma_H^+)$ is transitive, we prove that there is a dense set $D$ in $\mathbb X$ such that for each $s\in D$,  $\Cl(\mathcal U^{\oplus}_H(s))=\mathbb X$. Let $D=(0,1)\cup(2,3)\cup(4,5)\cup(6,7)\cup\ldots$ Then $D$ is dense in $\mathbb X$. Let $s\in D$ be any point and let $\ell$ be a non-negative integer such that $s\in (2\ell,2\ell+1)$. Note that 
	$$
	s,s-2,s-4, s-6,\ldots, s-2\ell\in \mathcal U^{\oplus}_H(s)
	$$
	and let $t=s-2\ell$. Then $t\in (0,1)$. It follows from the definition of $H$ that for all non-negative integers $m$, $n$ and $k$, 
	$$
	t^{\frac{2^m}{3^n}}+k\cdot 2\in \mathcal U^{\oplus}_H(t).
	$$
	It follows from Theorem \cite[Lemma 4.9, page 19]{BE} that $\big\{t^{\frac{2^m}{3^n}}+k\cdot 2 \ | \ m,,n,k\in \mathbb Z\big\}$ is dense in $\mathbb X$. Since 
	$$
	\Big\{t^{\frac{2^m}{3^n}}+k\cdot 2 \ | \ m,,n,k\in \mathbb Z\Big\}\subseteq \mathcal U^{\oplus}_H(t)\subseteq \mathcal U^{\oplus}_H(s),
	$$
	it follows that $\mathcal U^{\oplus}_H(s)$ is dense in $\mathbb X$. Therefore, by Theorem \ref{A}, $(\mathbb X_H^+,\sigma_H^+)$ is transitive
	\end{proof}
Next, we examine the space $\mathbb X_H$. 
\begin{definition}
	For each positive integer $k$, we use $L_k$ to denote
	$$
	L_k=\{(\ldots, t_{-2},t_{-1},t_0;t_1,t_2,\ldots)\in \mathbb X_H \ | \ t_0\in I_k\}.
	$$
\end{definition}
\begin{observation}
	Note that for each positive integer $k$, $L_k$ is compact and that
	$$
	\mathbb X_H=\left(\bigcup_{k=1}^{\infty}L_k\right)\cup \{(\ldots,\infty,\infty;\infty,\ldots)\}
	$$
\end{observation}
\begin{observation}
Let $k$ be a positive integer. Then for all $\mathbf x,\mathbf y\in L_k$,
$$
\textup{d}_{\mathbb X}(\mathbf x(0),\mathbf y(0))\leq \diam(I_k)=\frac{1}{2^{2k-1}}.
$$
Also, for all $\mathbf x,\mathbf y\in L_k$ and for each integer $n$,
$$
|n|< {2k-1} ~~~  \Longrightarrow  ~~~   \textup{d}_{\mathbb X}(\mathbf x(n),\mathbf y(n))\leq \diam(I_{k-n})=\frac{1}{2^{{2k-1}-|n|}}
$$
and
$$
|n|\geq {2k-1} ~~~  \Longrightarrow  ~~~   \textup{d}_{\mathbb X}(\mathbf x(n),\mathbf y(n))\leq \frac{1}{2}
$$
(in fact, the last inequality is true for each integer $n$). Therefore, for all $\mathbf x,\mathbf y\in L_k$,
\begin{align*}
	\frac{\textup{d}_{\mathbb X}(\mathbf x(0),\mathbf y(0))}{2^{|0|}}&=\textup{d}_{\mathbb X}(\mathbf x(0),\mathbf y(0))\leq \frac{1}{2^{2k-1}}.
\end{align*}
Also, for all $\mathbf x,\mathbf y\in L_k$ and for each integer $n$, if $|n|< {2k-1}$, then 
\begin{align*}
	\frac{\textup{d}_{\mathbb X}(\mathbf x(n),\mathbf y(n))}{2^{|n|}}&=\textup{d}_{\mathbb X}(\mathbf x(n),\mathbf y(n))\cdot \frac{1}{2^{|n|}}\leq \frac{1}{2^{k-|n|}}\cdot \frac{1}{2^{|n|}}=\frac{1}{2^{2k-1}},
	\end{align*}
and, if $|n|\geq  {2k-1}$, then 
\begin{align*}
\frac{\textup{d}_{\mathbb X}(\mathbf x(n),\mathbf y(n))}{2^{|n|}}&=\textup{d}_{\mathbb X}(\mathbf x(n),\mathbf y(n))\cdot \frac{1}{2^{|n|}}\leq \frac{1}{2}\cdot \frac{1}{2^{|n|}}<\frac{1}{2^{2k-1}}.
	\\
\end{align*}
It follows that for all $\mathbf x,\mathbf y\in L_k$ and for each integer $n$,
$$
\frac{\textup{d}_{\mathbb X}(\mathbf x(n),\mathbf y(n))}{2^{|n|}}\leq \frac{1}{2^{2k-1}}.
$$
Therefore, 
\begin{align*}
	&\diam(L_k)=\sup \{\textup{D}_{\mathbb X}(\mathbf x,\mathbf y) \ | \ \mathbf x,\mathbf y\in L_k\}=\\
	&\sup \Big\{\max \Big\{\frac{\textup{d}_{\mathbb X}(\mathbf x(n),\mathbf y(n))}{2^{|n|}} \ \big| \ n \textup{ is an integer}\Big\} \ | \ \mathbf x,\mathbf y\in L_k\Big\}\leq \\
	&\sup \Big\{\max \Big\{\frac{1}{2^{{2k-1}}} \ \big| \ k \textup{ is an integer}\Big\} \ | \ \mathbf x,\mathbf y\in L_k\Big\}=\frac{1}{2^{2k-1}}.
\end{align*}	
\end{observation}
\begin{definition}
We define the functions $f_{1,2},f_{1,3}:I_1\rightarrow \mathbb X$ as follows. For each $t\in I_1$, we define
$$
f_{1,2}(t)=t^{\frac{1}{3}}  ~~~  \textup{ and } ~~~  f_{1,3}(t)=t+2.
$$
We also define the functions $f_{2,1},f_{2,2},f_{2,3}:I_2\rightarrow \mathbb X$ as follows. For each $t\in I_2$, we define
$$
f_{2,1}(t)=t-2, ~~~ f_{2,2}(t)=t^2  ~~~  \textup{ and } ~~~  f_{2,3}(t)=t+2.
$$
Also, for each positive integer $k$, we define the functions $f_{k,1},f_{k,2},f_{k,3}:I_k\rightarrow \mathbb X$ as follows. For each $t\in I_k$, we define
$$
f_{k,1}(t)=t-2, ~~~ f_{k,2}(t)=t  ~~~  \textup{ and } ~~~  f_{k,3}(t)=t+2.
$$
We also use $\mathcal H$ to denote $\mathcal H=\{f_{1,2},f_{1,3}\}\cup \bigcup_{k=2}^{\infty}\{f_{k,1},f_{k,2},f_{k,3}\}$ (see Figure \ref{uncun} -- the relation $H$ contains as a subset the union of the graphs of the defined functions).
\end{definition}
\begin{observation}
	Let $\mathbf x\in \mathbb X_H\setminus \{(\ldots,\infty,\infty;\infty,\ldots)\}$. Then there is a unique point 	$\mathbf h=(\ldots,h_{-2},h_{-1},h_{0};h_1,h_2,\ldots)\in \prod_{k=-\infty}^{\infty}\mathcal H$ such that for each  integer $k$,
	$$
	\mathbf x(k+1)=h_k(\mathbf x(k)).
	$$
	Let $k$ be any integer. Then for each positive integer $\ell$ the following hold. 
	If $\ell=1$, then 
	\begin{enumerate}
		\item if $h_k=f_{\ell,2}$, then $h_{k+1}\in \{f_{\ell,2},f_{\ell,3}\}$. 
		\item if $h_k=f_{\ell,3}$, then $h_{k+1}\in \{f_{\ell+1,1},f_{\ell+1,2},f_{\ell+1,3}\}$. 
	\end{enumerate} 
	If $\ell=2$, then 
	\begin{enumerate}
		\item if $h_k=f_{\ell,1}$, then $h_{k+1}\in \{f_{\ell-1,2},f_{\ell-1,3}\}$. 
		\item if $h_k=f_{\ell,2}$, then $h_{k+1}\in \{f_{\ell,1},f_{\ell,2},f_{\ell,3}\}$. 
		\item if $h_k=f_{\ell,3}$, then $h_{k+1}\in \{f_{\ell+1,1},f_{\ell+1,2},f_{\ell+1,3}\}$. 
	\end{enumerate} 
		If $\ell>2$, then 
	\begin{enumerate}
		\item if $h_k=f_{\ell,1}$, then $h_{k+1}\in \{f_{\ell-1,1},f_{\ell-1,2},f_{\ell-1,3}\}$. 
		\item if $h_k=f_{\ell,2}$, then $h_{k+1}\in \{f_{\ell,1},f_{\ell,2},f_{\ell,3}\}$. 
		\item if $h_k=f_{\ell,3}$, then $h_{k+1}\in \{f_{\ell+1,1},f_{\ell+1,2},f_{\ell+1,3}\}$. 
	\end{enumerate} 
\end{observation}
\begin{definition}
	We define $\mathbf K$ to be the subset of the set $\prod_{k=-\infty}^{\infty}\mathcal H$, defined as follows. For any point $\mathbf h=(\ldots,h_{-2},h_{-1},h_{0};h_1,h_2,\ldots)\in \prod_{k=-\infty}^{\infty}\mathcal H$, $\mathbf h\in \mathbf K$ if and only if for each  integer $k$ and for each positive integer $\ell$ the following hold.
	\begin{enumerate}
		\item If $\ell=1$, then 
	\begin{enumerate}
		\item if $h_k=f_{\ell,2}$, then $h_{k+1}\in \{f_{\ell,2},f_{\ell,3}\}$. 
		\item if $h_k=f_{\ell,3}$, then $h_{k+1}\in \{f_{\ell+1,1},f_{\ell+1,2},f_{\ell+1,3}\}$.
	\end{enumerate} 
	\item If $\ell=2$, then 
	\begin{enumerate}
		\item if $h_k=f_{\ell,1}$, then $h_{k+1}\in \{f_{\ell-1,2},f_{\ell-1,3}\}$. 
		\item if $h_k=f_{\ell,2}$, then $h_{k+1}\in \{f_{\ell,1},f_{\ell,2},f_{\ell,3}\}$. 
		\item if $h_k=f_{\ell,3}$, then $h_{k+1}\in \{f_{\ell+1,1},f_{\ell+1,2},f_{\ell+1,3}\}$. 
	\end{enumerate} 
		\item If $\ell>2$, then 
	\begin{enumerate}
		\item if $h_k=f_{\ell,1}$, then $h_{k+1}\in \{f_{\ell-1,1},f_{\ell-1,2},f_{\ell-1,3}\}$. 
		\item if $h_k=f_{\ell,2}$, then $h_{k+1}\in \{f_{\ell,1},f_{\ell,2},f_{\ell,3}\}$. 
		\item if $h_k=f_{\ell,3}$, then $h_{k+1}\in \{f_{\ell+1,1},f_{\ell+1,2},f_{\ell+1,3}\}$. 
	\end{enumerate} 
	\end{enumerate}	
	We will also use $\mathbf h(j)$ to denote $\mathbf h(j)=h_j$.
\end{definition}
\begin{definition}
	For each positive integer $k$, we define 
	$$
	\mathbf K_k=\{\mathbf h\in \mathbf K \ | \ \mathbf h(0):I_k\rightarrow \mathbb X\}.
	$$
\end{definition}
\begin{observation}
	Note that 
	$$
	\mathbf K=\bigcup_{k=1}^{\infty}\mathbf K_k.
	$$
\end{observation}
\begin{observation}
	For each $\mathbf x\in \mathbb X_H$ and for each positive integer $k$, $\mathbf x\in L_k$ if and only if there is a unique point $\mathbf h=(\ldots,h_{-2},h_{-1},h_{0};h_1,h_2,\ldots)\in \mathbf K_k$ such that for each  integer $j$,
	$$
	\mathbf x(j+1)=h_j(\mathbf x(j)).
	$$
\end{observation}

\begin{definition}
Let $f_{1,1}=f_{1,2}$, let $k$ be a positive integer and let for each integer $\ell$, 
$$
A_{\ell}=\bigcup_{i=\max\{k-|\ell|,1\}}^{k+|\ell|}\{f_{i,1},f_{i,2},f_{i,3}\}.
$$ 
We equip each $A_\ell$ with the discrete topology. Then we use $\mathbf C_k$ to denote the set 
	\begin{align*}
	 \mathbf C_k=\prod_{\ell=-\infty}^{\infty}A_{\ell}.
	\end{align*}
\end{definition}
\begin{observation}
	For each positive integer $k$, $\mathbf C_k$ is a Cantor set  and $\mathbf K_k\subseteq \mathbf C_k$.  
\end{observation}
\begin{definition}
	For each positive integer $k$ and for each $\mathbf h\in \mathbf K_k$, we define
	$$
	A_{k,\mathbf h}=\{\mathbf x\in L_k \ | \ \textup{for each integer j}, \mathbf x(j+1)=\mathbf h(j)(\mathbf x(j))\}.
	$$
\end{definition}
\begin{observation}\label{lik}
Note that for each positive integer $k$,
$$
L_k=\bigcup_{\mathbf h\in \mathbf K_k}A_{k,\mathbf h}.
$$
	Let $k$ be a positive integer and let $\mathbf h\in \mathbf K_k$. Since for each integer $j$, $\mathbf h(j)$ is an increasing  homeomorphism from an interval $I_m$ to an interval $I_n$, it follows that $A_{k,\mathbf h}$ is an arc in $L_k$. For each of the end-points $\mathbf e$ of the arc $A_{k,\mathbf h}$, either all coordinates of $\mathbf e$ are odd or all coordinates of $\mathbf e$ are even. Also, note that  for all $\mathbf h_1,\mathbf h_2\in \mathbf K_k$,
	$$
	\mathbf h_1\neq \mathbf h_2 ~~~  \Longrightarrow  ~~~   A_{k,\mathbf h_1}\cap A_{k,\mathbf h_2}=\emptyset.
	$$  
\end{observation}

\begin{theorem}\label{muca}
	Let $k$ be a positive integer. Then  $\mathbf K_k$ is a closed subset of $\mathbf C_k$ and $L_k$ is homeomorphic to $\mathbf K_k\times [0,\frac{1}{2^{2k-1}}]$. 
\end{theorem}
\begin{proof}
	First, we show that $L_k$ is homeomorphic to $\mathbf K_k\times [0,\frac{1}{2^{2k-1}}]$. Let $\varphi:\mathbf K_k\times [0,\frac{1}{2^{2k-1}}]\rightarrow L_k$ be defined by
$$
\varphi(\mathbf h,t)=(\ldots,t_{-2},t_{-1},t_0;t_1,t_2,\ldots)
$$
for any $(\mathbf h,t)\in \mathbf K_k\times [0,\frac{1}{2^{2k-1}}]$, where $t_0=2^{2k-1}\cdot t+2k-2$ and for each integer $j$, $$t_{j+1}=\mathbf h(j)(t_j).$$ 
Then $\varphi$ is a homeomorphism. Since $L_k$ is compact, it follows that $\mathbf K_k\times [0,\frac{1}{2^{2k-1}}]$ is compact. Therefore, $\mathbf K_k$ is a closed subset of $\mathbf C_k$.
\end{proof}
\begin{observation}\label{jure}
	Note that in the proof of Theorem \ref{muca}, the homeomorphism $\varphi$ is constructed in such a way that for each $\mathbf x\in L_k$, if all the coordinates of $\mathbf x$ are even, then $p_2(\varphi^{-1}(\mathbf x))=0$. 
\end{observation}
\begin{theorem}
	For each positive integer $k$, $\mathbf K_k$ is a Cantor set. 
\end{theorem}
	\begin{proof}
		Suppose that there is a positive integer $k$ such that $\mathbf K_k$ is not a Cantor set. Note that $\mathbf K_k$ is a totally disconnected metric compactum since by Theorem \ref{muca} it is a closed subset of a Cantor set.  Since $\mathbf K_k$ is not a Cantor set, it follows from  \cite[Theorem 7.14, page 109]{nadler} there is an isolated point in $\mathbf K_k$. Let $\mathbf h\in \mathbf K_k$ be an isolated point of $\mathbf K_k$. Then $\mathbf h$ is an isolated point in $\mathbf K$. It follows that $A_{k,\mathbf h}$ is an isolated arc in $\mathbb X_H$ (meaning that there is an open set $U$ in $\mathbb X_H$ such that $A_{k,\mathbf h}\subseteq U$ and $(\mathbb X_H\setminus A_{k,\mathbf h})\cap U=\emptyset$). 
		
		Let $\mathbf x\in \mathbb X_H$ be any transitive point in $(\mathbb X_H,\sigma_H)$. If $\mathbf x$ is an element of an isolated arc $A$ in $\mathbb X_H$, then (since by Theorem \ref{lejgaungatamlehre}, $\sigma_H$ is a transitive homeomorphism) there is a positive integer $n$ such that $\sigma_H^n(\mathbf x)\in A$ and, therefore, $\sigma_H^n(A)=A$. It follows that $\mathbb X_H$ is the union of $n$ mutually disjoint arcs: $\sigma_H(A)$, $\sigma_H^2(A)$, $\sigma_H^3(A)$, $\ldots$, $\sigma_H^n(A)$, which is a contradiction. It follows that $\mathbf x$ is not an element of an isolated arc. Let $A$ be an isolated arc in $\mathbb X_H$ and let $U$ be an open set in $\mathbb X_H$ such that $A\subseteq U$ and $(\mathbb X_H\setminus A)\cap U=\emptyset$. Then for each non-negative integer $k$, $\sigma_H^k(\mathbf x)\not \in U$. It follows that $\mathbf x$ is not a transitive point in $(\mathbb X_H,\sigma_H)$, which is a contradiction. 
	\end{proof}


\begin{definition}
	Let $C$ be the standard middle-third Cantor set in $[0,1]$. For each positive integer $k$, we use $C_k$ to denote 
	$C_k=C\cap [c_k,d_k]$, 	where $c_1=0$, $d_1=\frac{1}{3}$, and for each positive integer $k$, $c_{k+1}=d_{k}+\frac{1}{3^k}$ and $d_{k+1}=c_{k+1}+\frac{1}{3^{k+1}}$. 	  
\end{definition}
\begin{observation}
	Note that for each positive integer $k$, $C_k$ is a Cantor set and that 
	$$
		C=\left(\bigcup_{k=1}^{\infty}C_k\right)\cup \{1\}.
		$$
		Also, note that for all positive integers $k$ and $\ell$,
		$$
		k\neq \ell ~~~  \Longrightarrow  ~~~  C_k\cap C_{\ell}=\emptyset.
		$$
\end{observation}
In the following theorem, we obtain a model for our two-sided Mahavier product $\mathbb X_H$. This will be used later in Theorem \ref{dud} where we show that the members of our uncountable family are in fact smooth fans. 
\begin{theorem}\label{ulci}
There is a homeomorphism 
$$
\varphi:\mathbb X_H\rightarrow \left(\bigcup_{k=1}^{\infty}\Big(C_k\times \Big[0,\frac{1}{2^{2k-1}}\Big]\Big)\right)\cup \{(1,0)\}
$$ 
such that for each $\mathbf x\in \mathbb X_H$, if all the coordinates of $\mathbf x$ are even, then $\varphi(\mathbf x)=(c,0)$ for some $c\in C$.
\end{theorem}
\begin{proof}
	For each positive integer $k$, let 
	$$
	f_k:\mathbf K_k\rightarrow C_k
	$$
	 be a homeomorphism. By Theorem \ref{muca}, each $L_k$ is homeomorphic to $\mathbf K_k\times [0,\frac{1}{2^{2k-1}}]$. For each positive integer $k$, let 
	   $$
	   v_k:L_k\rightarrow \mathbf K_k\times \Big[0,\frac{1}{2^{2k-1}}\Big]
	   $$
	    be a homeomorphism such that for each $\mathbf x\in L_k$, if all the coordinates of $\mathbf x$ are even, then $p_2(v_k(\mathbf x))=0$ (such a homeomorphism does exist by Observation \ref{jure}). 
Then 
$$
\varphi:\mathbb X_H\rightarrow \left(\bigcup_{k=1}^{\infty}\Big(C_k\times \Big[0,\frac{1}{2^{2k-1}}\Big]\Big)\right)\cup \{(1,0)\},
$$
 defined by 
$$
\varphi(\mathbf x)=\begin{cases}
				(1,0)\text{;} & \mathbf x=(\ldots, \infty,\infty;\infty,\ldots) \\
				(f_k(p_1(v_k(\mathbf x))),p_2(v_k(\mathbf x)))\text{;} & \textup{there is a positive integer } k \textup{ such that } \mathbf x\in L_k.
			\end{cases}
$$
for each $\mathbf x\in \mathbb X_H$, is a homeomorphism such that for each $\mathbf x\in \mathbb X_H$, if all the coordinates of $\mathbf x$ are even, then $\varphi(\mathbf x)=(c,0)$ for some $c\in C$. See Figure \ref{mahi}, where the space $\left(\bigcup_{k=1}^{\infty}\Big(C_k\times \Big[0,\frac{1}{2^{2k-1}}\Big]\Big)\right)\cup \{(1,0)\}$ is presented. 
\begin{figure}[h!]
	\centering
		\includegraphics[width=25em]{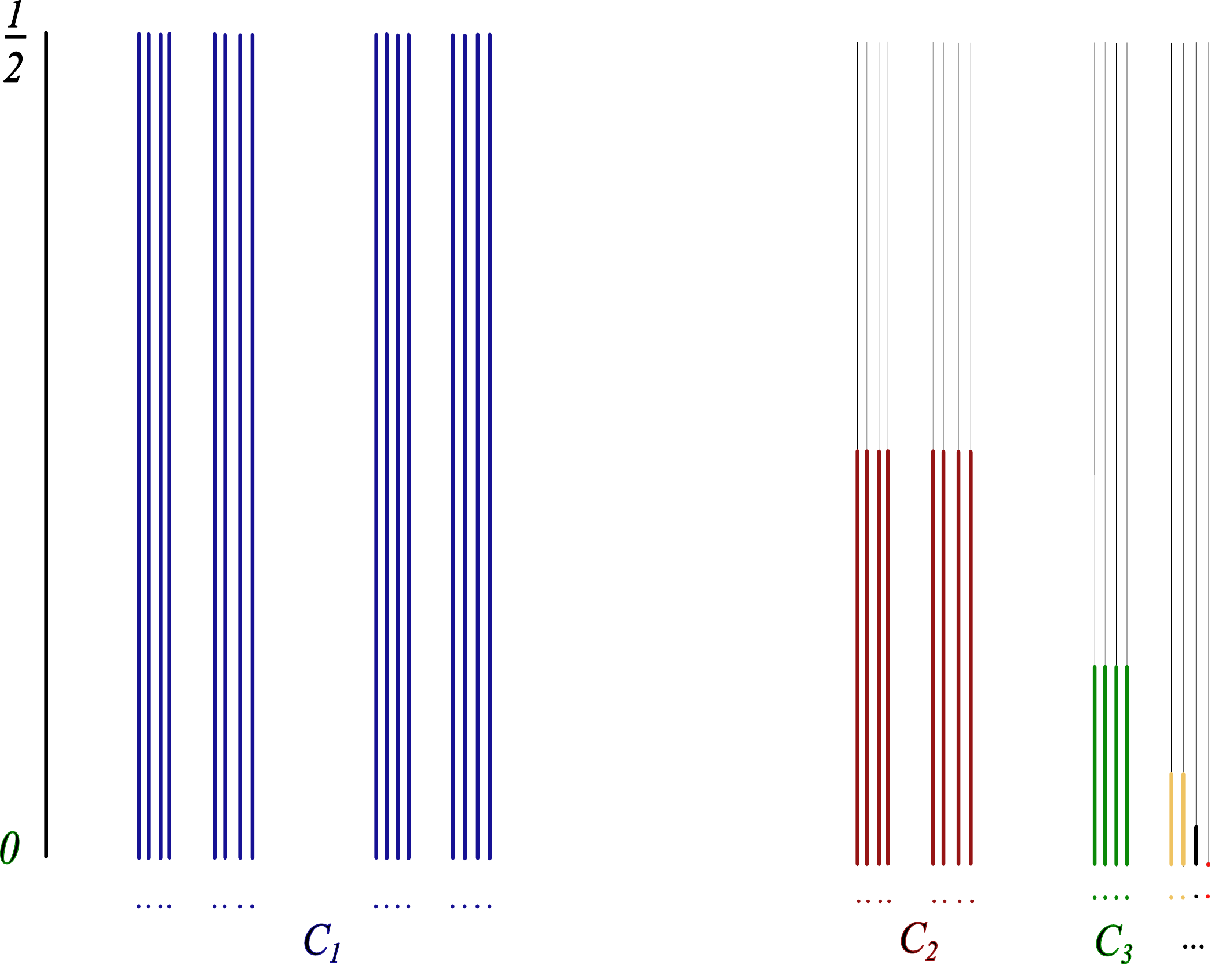}
	\caption{The space $\displaystyle \left(\bigcup_{k=1}^{\infty}\Big(C_k\times \Big[0,\frac{1}{2^{2k-1}}\Big]\Big)\right)\cup \{(1,0)\}$}
	\label{mahi}
\end{figure}
\end{proof}
\begin{definition}
	We choose and fix one of the homeomorphisms 
	$$
\varphi:\mathbb X_H\rightarrow \left(\bigcup_{k=1}^{\infty}\Big(C_k\times \Big[0,\frac{1}{2^{2k-1}}\Big]\Big)\right)\cup \{(1,0)\}
$$ 
such that for each $\mathbf x\in \mathbb X_H$, if all the coordinates of $\mathbf x$ are even, then $\varphi(\mathbf x)=(c,0)$ for some $c\in C$, and we denote it by $\varphi_0$.
\end{definition}

Next, we use the space $\mathbb X_H$ and the model that we obtained in Theorem \ref{ulci} to construct a family of uncountable many pairwise non-homeomorphic smooth fans. First, we introduce the following definitions.

\begin{definition}
We use $\mathbb A$ to denote the product 
$$
\mathbb A=\{1,2\}\times \{3,4\}\times\{5,6\}\times\{7,8\}\times\{9,10\}\times\ldots
$$
\end{definition}
\begin{observation}\label{uncunc}
	Note that $\mathbb A$ is uncountable.
\end{observation}
 Using the set $\mathbb A$, we define three relations on $\mathbb X_H$; see Definitions \ref{Hum1}, \ref{Hum2} and \ref{destil}.
\begin{definition}\label{Hum1}
	We define the relation $\approx$ on $\mathbb X_H$ as follows: for all $\mathbf x,\mathbf y\in \mathbb X_H$, we define $\mathbf x\approx \mathbf y$ if and only if one of the following holds:
	\begin{enumerate}
	\item\label{11} $\mathbf x=\mathbf y$,
	\item\label{21} $p_2(\varphi_0(\mathbf x))=0$ and $\varphi_0(\mathbf y)=(1,0)$ or and $\varphi_0(\mathbf x)=(1,0)$ and $p_2(\varphi_0(\mathbf y))=0$, 
	\item\label{31} $p_2(\varphi_0(\mathbf x))=p_2(\varphi_0(\mathbf y))=0$; 
	\end{enumerate}
	 see Figure \ref{urh1}.
	\begin{figure}[h!]
\centering
	\includegraphics[width=30em]{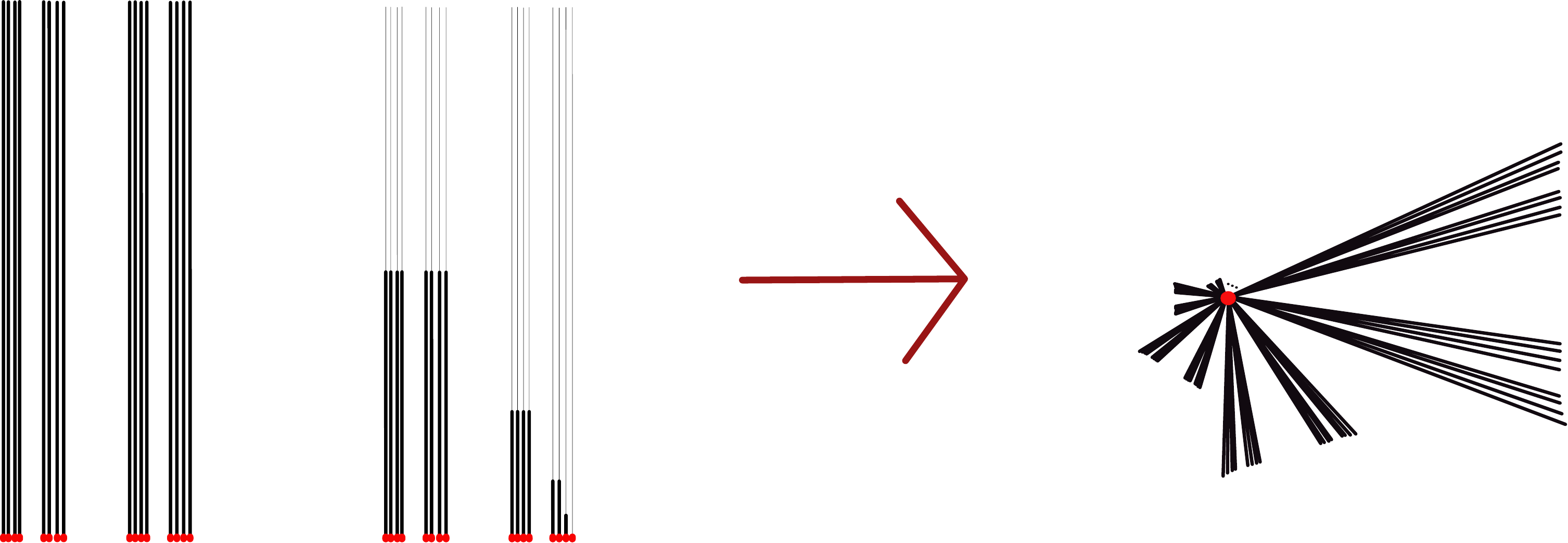}
\caption{The relation $\approx $ from Definition \ref{Hum1}}
\label{urh1}
\end{figure} 
\end{definition}
\begin{definition}
For each positive integer $k$, we use $M_k$ to be the following subspace of $L_k$:
$$
M_k=\{(\ldots,t,t;t,\ldots) \ | \ t\in I_k\}.
$$
\end{definition}
Note that there are infinitely many of these $M_k$'s, and that plays a necessary role in creating the uncountable family of smooth fans. Also, note that the shift map restricted to these $M_k$'s is the identity, and that is crucial in maintaining the transitivity of the induced map after identifications using the following relation.
\begin{definition}\label{Hum2}
	Let $\mathbf a= (a_1,a_2,a_3,\ldots)\in \mathbb A$. Then we define the relation $\approx_{\mathbf a}$ on $\mathbb X_H$ as follows: for all $\mathbf x,\mathbf y\in \mathbb X_H$, we define $\mathbf x\approx_{\mathbf a} \mathbf y$ if and only if one of the following holds:
	\begin{enumerate}
	\item\label{12} $\mathbf x=\mathbf y$,
	\item \label{22} there is a positive integer $k$ and there is an $i\in \{1,2,3,\ldots,a_k\}$ such that either
	\begin{enumerate}
	\item  $\mathbf x\in M_{k^2+2}$ and $\mathbf y\in M_{k^2+2+i}$, and
	\item  $p_2(\varphi_0(\mathbf x))=p_2(\varphi_0(\mathbf y))$
	\end{enumerate}
	or 
	\begin{enumerate}
	\item  $\mathbf y\in M_{k^2+2}$ and $\mathbf x\in M_{k^2+2+i}$, and
	\item  $p_2(\varphi_0(\mathbf x))=p_2(\varphi_0(\mathbf y))$;
	\end{enumerate}
	\end{enumerate}
	see Figure 
	\ref{JuMa}, which illustrates how arcs $M_{k^2+2+1}$, $M_{k^2+2+2}$, $M_{k^2+2+3}$, $\ldots$, $M_{k^2+2+a_k}$ are being glued to the arc $M_{k^2+2}$.
\end{definition}
\begin{figure}[h!]
\centering
	\includegraphics[width=32em]{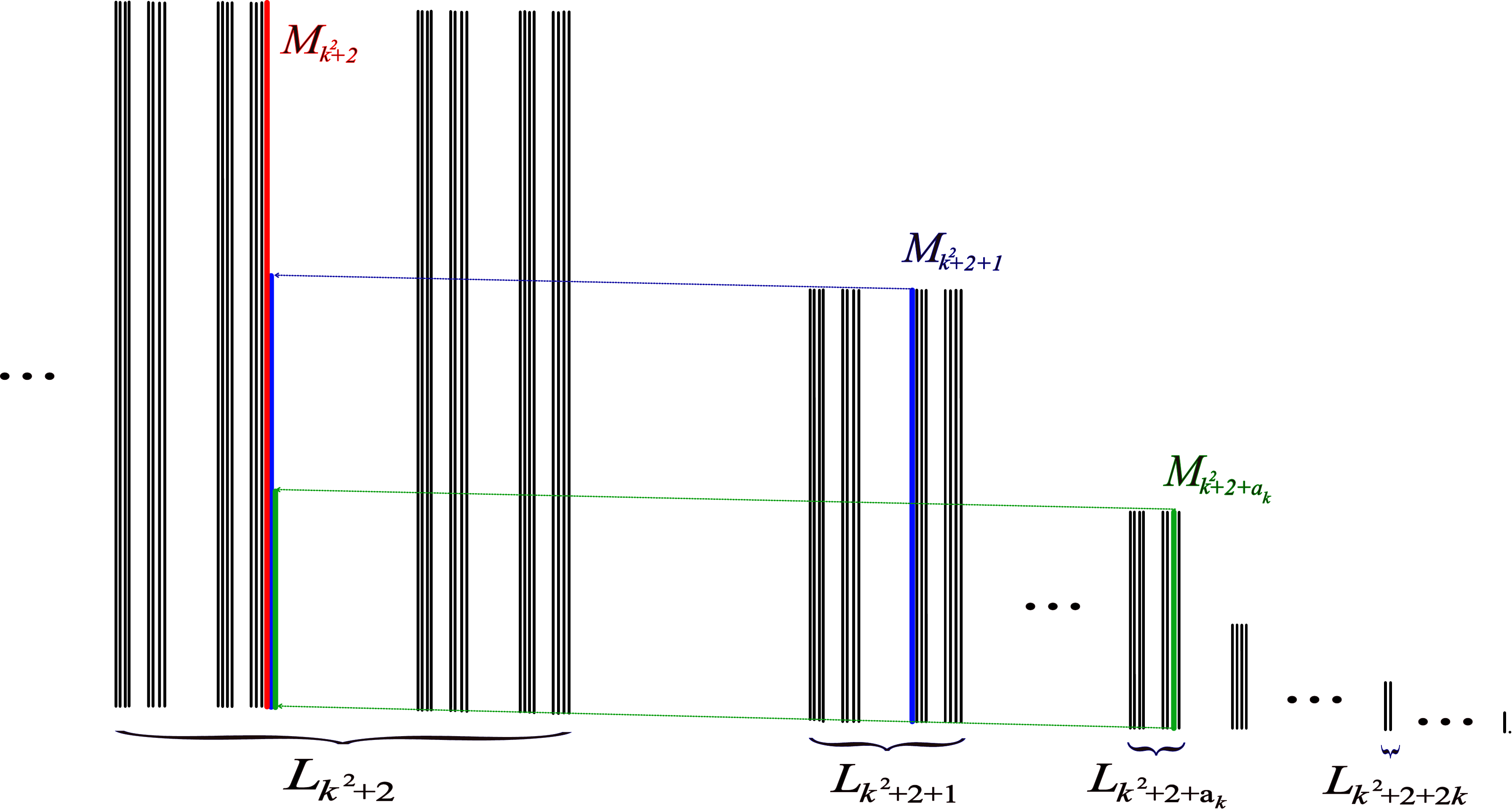}
\caption{The relation $\approx_{\mathbf a}$ from Definition \ref{Hum1}}
\label{JuMa}
\end{figure}

\begin{definition}\label{destil}
	For each $\mathbf a= (a_1,a_2,a_3,\ldots)\in \mathbb A$,  we define the relation $\sim_{\mathbf a}$ on $\mathbb X_H$ by
			$$
	\mathbf x\sim_{\mathbf a} \mathbf y ~~~  \Longleftrightarrow  ~~~  \mathbf x\approx \mathbf y \textup{ or there is } \mathbf a\in \mathbb A \textup{ such that } \mathbf x\approx_{\mathbf a} \mathbf y
	$$
	 for all $\mathbf x,\mathbf y\in \mathbb X_H$.
\end{definition}
\begin{observation}
	Note that $\sim_{\mathbf a}$ is an equivalence relation on $\mathbb X_H$.
\end{observation}
\begin{definition}
	For each $\mathbf a\in \mathbb A$, we use $F_{\mathbf a}$ to denote  the quotient space
	$$
	F_{\mathbf a}=\mathbb X_H/_{\sim_{\mathbf a}}.
	$$
\end{definition}

\begin{theorem}\label{dud}
For each $\mathbf a\in \mathbb A$,  $F_{\mathbf a}$ is a smooth fan. 
\end{theorem}
\begin{proof}
First, let $\mathbf a=(a_1,a_2,a_3,\ldots)\in \mathbb A$ and let 
  $$
  i:\left(\bigcup_{k=1}^{\infty}\Big(C_k\times \Big[0,\frac{1}{2^{2k-1}}\Big]\Big)\right)\cup \{(1,0)\}\rightarrow C\times [0,1]
  $$
   be the inclusion function. For each positive integer $k$ and for each integer $i\in \{0,1,2,3,\ldots,a_k\}$, let
$A_{k,i}$ be the connected component of $C\times [0,1]$ such that 
$$
i(\varphi_0(M_{k^2+2+i}))\subseteq A_{k,i};
$$
 see Figure \ref{spodaj}.

\begin{figure}[h!]
\centering
	\includegraphics[width=35em]{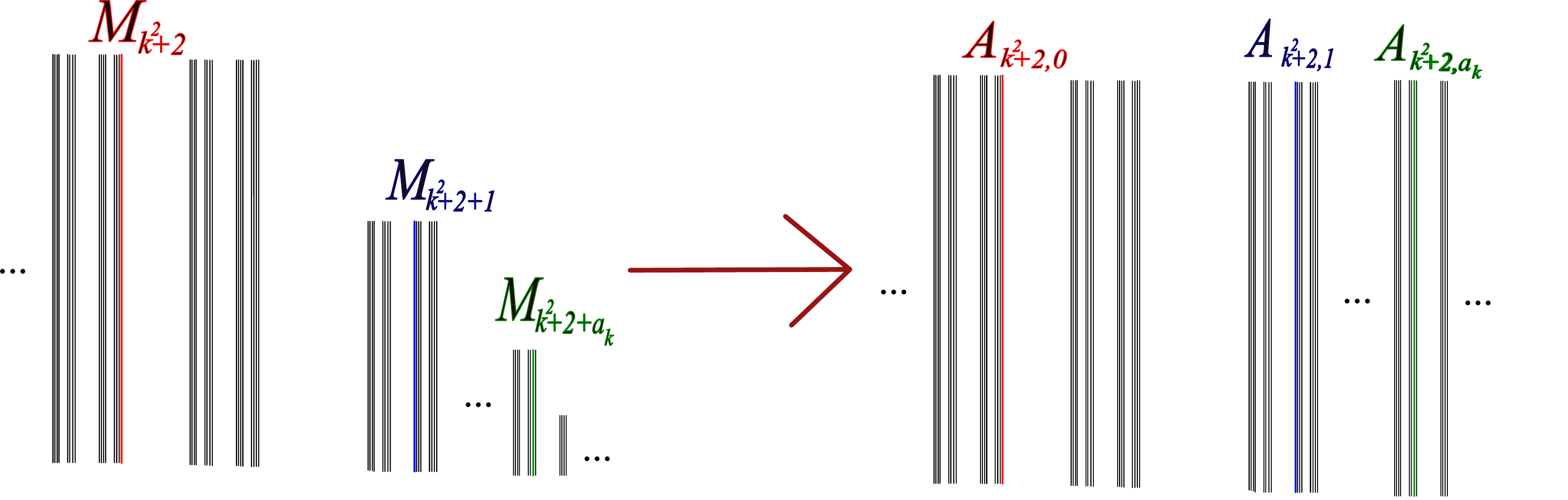}
\caption{The arcs $A_{k,i}$ in $C\times [0,1]$}
\label{spodaj}
\end{figure} 
Next, let $\sim_1$ be the equivalence relation on $C\times [0,1]$, defined as follows. For all $(c_1,t_1),(c_2,t_2)\in C\times [0,1]$, we define that $(c_1,t_1)\sim_1 (c_2,t_2)$ if and only if one of the following holds. 
\begin{enumerate}
	\item $(c_1,t_1)=(c_2,t_2)$,
	\item $t_1=t_2=0$.
\end{enumerate}
Note that $(C\times [0,1])/_{\sim_1}$ is a Cantor fan. We also define $\sim_2$ to be the equivalence relation on $C\times [0,1]$, defined as follows. For all $(c_1,t_1),(c_2,t_2)\in C\times [0,1]$, we define that $(c_1,t_1)\sim_2 (c_2,t_2)$ if and only if one of the following holds. 
\begin{enumerate}
	\item $(c_1,t_1)=(c_2,t_2)$,
	\item there is a positive integer $k$ and an integer $i\in \{1,2,3,\ldots,a_k\}$ such that either 
	\begin{enumerate}
	\item $(c_1,t_1)\in A_{k^2+2,0}$ and $(c_2,t_2)\in A_{k^2+2+i}$, and
	\item $t_1=t_2$,
	\end{enumerate}
	or 
	\begin{enumerate}
	\item $(c_2,t_2)\in A_{k^2+2,0}$ and $(c_1,t_1)\in A_{k^2+2+i}$, and
	\item $t_1=t_2$.
	\end{enumerate}
\end{enumerate}
Finally, we define the equivalence relation $\sim$ on $C\times [0,1]$ as follows. For all $(c_1,t_1),(c_2,t_2)\in C\times [0,1]$, we define 
$$
(c_1,t_1)\sim (c_2,t_2) ~~~ \Longleftrightarrow ~~~  (c_1,t_1)\sim_1 (c_2,t_2) \textup{ or } (c_1,t_1)\sim_2 (c_2,t_2).
$$
Next, let 
$$
r:C\times [0,1]\rightarrow (C\times [0,1])/_{\sim_1}
$$
 be the quotient map defined by  
$$
r(c,t)=[(c,t)]_{\sim_1}=\{(d,s)\in C\times [0,1] \ | \  (d,s)\sim_1 (c,t)\}
$$
 for each $(c,t)\in C\times [0,1]$ and let
$$
q:C\times [0,1]\rightarrow (C\times [0,1])/_{\sim}
$$ 
be the quotient map defined by  
$$
q(c,t)=[(c,t)]_{\sim}=\{(d,s)\in C\times[0,1] \ | \  (d,s)\sim (c,t)\}
$$
 for each $(c,t)\in C\times [0,1]$. We use $F$ to denote $F=(C\times [0,1])/_{\sim}$.  Let 
 $$
 g:(C\times [0,1])/_{\sim_1}\rightarrow F
 $$
 be defined by
 $$
 g([(c,t)]_{\sim_1})=q(r^{-1}([(c,t)]_{\sim_1}))
 $$
 for any $(c,t)\in C\times [0,1]$. Note that $g$ is a well-defined confluent surjection.  Since $(C\times [0,1])/_{\sim_1}$ is a smooth fan (in fact, it is a Cantor fan), it follows from Theorem \ref{char} that $F$ is also a smooth fan. 
Finally, let
 $$
 p:\mathbb X_H\rightarrow F_{\mathbf a}
 $$
  be the quotient map defined by  
$$
p(\mathbf x)=[\mathbf x]=\{\mathbf y\in \mathbb X_H \ | \  \mathbf y\sim_{\mathbf a} \mathbf x\}
$$
 for each $\mathbf x\in \mathbb X_H$. Then 
$$
f:F_{\mathbf a}\rightarrow F,
$$
defined by 
$$
f([\mathbf x])=q(i(\varphi_0(p^{-1}([\mathbf x]))))
$$
for any $\mathbf x\in \mathbb X_H$, is an embedding of $F_{\mathbf a}$ into the smooth fan $F$. Therefore, $F_{\mathbf a}$ is a smooth fan. 
\end{proof}
In Theorem \ref{gor}, we prove that each $F_{\mathbf a}$ admits a transitive homeomorphism. In it, we use the following observation.
	\begin{observation}\label{sito}
	Let $\mathbf a\in \mathbb A$.  For all $\mathbf x,\mathbf y\in \mathbb X_H$,
	$$
	\mathbf x\sim_{\mathbf a} \mathbf y  ~~~  \Longleftrightarrow  ~~~  \sigma_H(\mathbf x)\sim_{\mathbf a} \sigma_H(\mathbf y).
	$$
\end{observation}
\begin{theorem}\label{gor}
Let $\mathbf a\in \mathbb A$.  The mapping $\sigma_{H}^{\star}:F_{\mathbf a}\rightarrow F_{\mathbf a}$, defined by 
			$$
			\sigma_{H}^{\star}([\mathbf x])=[\sigma_H(\mathbf x)]
			$$
 for each $\mathbf x\in \mathbb X_H$, is a transitive homeomorphism. 
\end{theorem}
\begin{proof}
 By Theorem \ref{lejgaungatamlehre} and Observation \ref{juju}, $\sigma_H$ is a transitive homeomorphism. It follows from Lemma \ref{sito} and from Proposition \ref{kvocienti}, that $\sigma_{H}^{\star}$ is a transitive homeomorphism. 
 \end{proof}
 \begin{observation}
 	Note that for each positive integer $k$,  this transitive homeomorphism $\sigma_{H}^{\star}$, restricted to $M_k/_{\sim_\mathbf a}=\{[\mathbf x] \ | \ \mathbf x\in M_k\}$, is just the identity. 
 \end{observation}
\begin{definition}
	We use $\mathcal F$ to denote the family
	$$
	\mathcal F= \{F_{\mathbf a} \ | \ \mathbf a\in \mathbb A\}.
	$$
\end{definition}
By Theorems \ref{dud} and \ref{gor}, each member of $\mathcal F$ is a smooth fan that admits a transitive homeomorphism. Recall that by Obsrvation \ref{uncunc}, $\mathbb A$ is uncountable. So, if we show that for all $\mathbf a,\mathbf b\in \mathbb A$,  
$$
\mathbf a\neq \mathbf b ~~~ \Longrightarrow  ~~~  F_{\mathbf a} \textup{ and } F_{\mathbf b} \textup{ are not homeomorphic},
$$
then this proves that $\mathcal F$ is a family of uncountably many pairwise non-homeo\-morphic smooth fans  that admit transitive homeomorphisms. In the following definition, we define the new concept of JuMas, which will be used to prove this.
\begin{definition}
	Let $X$ be a fan with the top $o$. We define the set $\textup{JuMa}(X)$ as follows: 
	$$
	\textup{JuMa}(X)=\{x\in X\setminus \{o\} \ | \ \textup{ there is a sequence } (e_n) \textup{ in } E(X) \textup{ such that } \lim_{n\to \infty}e_n=x\}.  
$$
\end{definition} 
\begin{definition}
	Let $X$ be a fan with the top $o$. For each $e\in E(X)$, we use $A_X[o,e]$ to denote the arc in $X$ from $o$ to $e$. 
\end{definition}
\begin{proposition}\label{juma}
	Let $X$ and $Y$ be fans with tops $o_X$ and $o_Y$, respectively,  and let $f:X\rightarrow Y$ be a homeomorphism. Then for each $e\in E(X)$, 
	$$
	|A_X[o_X,e]\cap \textup{JuMa}(X)|=|A_Y[o_Y,f(e)]\cap\textup{JuMa}(Y)|.
	$$
	Here $|S|$ denotes the cardinality of $S$ for any set $S$.
\end{proposition}
\begin{proof}
	The lemma follows from the fact  that for each $x\in X$, 
	$$
	x\in \textup{JuMa}(X) ~~~ \Longrightarrow ~~~ f(x)\in \textup{JuMa}(Y),
	$$ 
	which is easy to see and we leave the details to the reader.
\end{proof}
\begin{corollary}\label{miodio}
	Let $X$ and $Y$ be fans with tops $o_X$ and $o_Y$, respectively. If there is $e\in E(X)$ such that for each $f\in E(Y)$,
	$$
	 |A_Y[o_Y,f]\cap\textup{JuMa}(Y)|\neq |A_X[o_X,e]\cap \textup{JuMa}(X)|,
	$$
	then $X$ and $Y$ are not homeomorphic. 
\end{corollary}
\begin{proof}
	The corollary follows directly from Proposition \ref{juma}.
\end{proof}
\begin{theorem}\label{sus}
For all $\mathbf a,\mathbf b\in \mathbb A$,  
$$
\mathbf a\neq \mathbf b ~~~ \Longrightarrow  ~~~  F_{\mathbf a} \textup{ and } F_{\mathbf b} \textup{ are not homeomorphic}.
$$
\end{theorem}
\begin{proof}
	Let $\mathbf a,\mathbf b\in \mathbb A$ be such that $\mathbf a\neq \mathbf b$. Let $\mathbf o_{\mathbf a}$ and $\mathbf o_{\mathbf b}$ be the tops of the fans $F_{\mathbf a}$ and $F_{\mathbf b}$, respectively. Since $\mathbf a\neq \mathbf b$,  there is a positive integer $k$ such that $\mathbf a(k)\neq \mathbf b(k)$. Then either $\mathbf a(k)=2k-1$ and $\mathbf b(k)=2k$ or $\mathbf a(k)=2k$ and $\mathbf b(k)=2k-1$. Without loss of generality we assume that  $\mathbf a(k)=2k-1$ and $\mathbf b(k)=2k$. It follows from the definition  of the relation $\sim_{\mathbf a}$ that in $F_{\mathbf a}$, there is an end-point $\mathbf e\in E(F_{\mathbf a})$  such that 
	$$
	|A_{F_{\mathbf a}}[\mathbf o_{\mathbf a},\mathbf e]\cap \textup{JuMa}(F_{\mathbf a})|=2k-1.
	$$
	Note that it follows from the definition  of the relation $\sim_{\mathbf b}$ that for each $\mathbf f\in E(F_{\mathbf b})$,  
	$$
	|A_{F_{\mathbf b}}[\mathbf o_{\mathbf b},\mathbf f]\cap \textup{JuMa}(F_{\mathbf b})|\neq 2k-1.
	$$
	Therefore, by Corollary \ref{miodio}, $F_{\mathbf a}$ and $F_{\mathbf b}$ are not homeomorphic.
\end{proof}
Finally, we state and prove Theorem \ref{mine} -- the main theorem of the paper.
\begin{theorem}\label{mine}
There is a family of uncountable many pairwise non-homeomorphic smooth fans that admit transitive homeomorphisms.
\end{theorem}
\begin{proof}
		By Theorem \ref{sus}, $\mathcal F$ is uncountable, since $\mathbb A$ is uncountable. By Theorem \ref{dud}, for each $\mathbf a\in \mathbb A$, $F_{\mathbf a}$ is a smooth fan and by Theorem \ref{gor},  $\sigma_H^*$ is a transitive homeomorphism on $F_{\mathbf a}$. This completes the proof.
\end{proof}
The following open problem is a good place to finish our paper.
\begin{problem}
	Is there a smooth fan $X$ with the top $o$ that has the following properties?
	\begin{enumerate}
		\item $X$ does not admit a transitive homeomorphism.
		\item For each $\varepsilon >0$, for each  $e\in E(X)$ and for each $x\in A_X[o,e]$, there is $f\in E(X)\setminus \{e\}$ such that 
		$$
		B(x,\varepsilon)\cap A_X[o,f]\neq \emptyset.
		$$  
	\end{enumerate}
\end{problem}
\section{Acknowledgement}
This work is supported in part by the Slovenian Research Agency (research projects J1-4632, BI-HR/23-24-011, BI-US/22-24-086 and BI-US/22-24-094, and research program P1-0285). 
	

\noindent I. Bani\v c\\
              (1) Faculty of Natural Sciences and Mathematics, University of Maribor, Koro\v{s}ka 160, SI-2000 Maribor,
   Slovenia; \\(2) Institute of Mathematics, Physics and Mechanics, Jadranska 19, SI-1000 Ljubljana, 
   Slovenia; \\(3) Andrej Maru\v si\v c Institute, University of Primorska, Muzejski trg 2, SI-6000 Koper,
   Slovenia\\
             {iztok.banic@um.si}           
     
				\-
				
		\noindent G.  Erceg\\
             Faculty of Science, University of Split, Rudera Bo\v skovi\' ca 33, Split,  Croatia\\
{{gorerc@pmfst.hr}       }    

                 	\-
					
  \noindent J.  Kennedy\\
             Department of Mathematics,  Lamar University, 200 Lucas Building, P.O. Box 10047, Beaumont, Texas 77710 USA\\
{{kennedy9905@gmail.com}       }    

	\-
				
		\noindent C.  Mouron\\
             Rhodes College,  2000 North Parkway, Memphis, Tennessee 38112  USA\\\
{{mouronc@rhodes.edu}       }    

                 	\-
				
		\noindent V.  Nall\\
             Department of Mathematics,  University of Richmond, Richmond, Virginia 23173 USA\\
{{vnall@richmond.edu}       }   



\end{document}